\documentclass[preprint,authoryear]{elsarticle}

\biboptions{sort&compress}

\usepackage{graphicx}      
\usepackage{natbib}        
\usepackage{amssymb}
\usepackage{amsmath}
\usepackage{enumerate}
\usepackage{placeins}
\usepackage{hyphenat}
\usepackage{epstopdf}
\usepackage{gensymb}

\newtheorem{thm}{Theorem}
\newtheorem{lem}{Lemma}

\newproof{pf}{Proof}

\begin{document}
\bibliographystyle{model5-names}
\pagestyle{plain}
\setcounter{page}{0}
\pagenumbering{arabic}
\begin{frontmatter}

\title{Sufficient Conditions for Feasibility and Optimality of Real-Time Optimization Schemes -- II. Implementation Issues} 

\author[epfl]{Gene A. Bunin}
\ead{gene.bunin@epfl.ch}
\author[epfl]{Gr\'egory Fran\c cois}
\ead{gregory.francois@epfl.ch}
\author[epfl]{Dominique Bonvin\corref{cor1}}
\ead{dominique.bonvin@epfl.ch}

\cortext[cor1] {Author to whom correspondence should be addressed (tel.: +41 21 6933843, fax: +41 21 6932574).}

\tnotetext[abb]{Nonstandard abbreviations used in this article: SCFO (sufficient conditions for feasibility and optimality).}

\address[epfl]{Laboratoire d'Automatique, Ecole Polytechnique F\'ed\'erale de Lausanne, CH-1015 Lausanne, Switzerland}

\begin{abstract}                

The idea of iterative process optimization based on collected output measurements, or ``real-time optimization'' (RTO), has gained much prominence in recent decades, with many RTO algorithms being proposed, researched, and developed. While the essential goal of these schemes is to drive the process to its true optimal conditions without violating any safety-critical, or ``hard'', constraints, no generalized, unified approach for \emph{guaranteeing} this behavior exists. In this two-part paper, we propose an implementable set of conditions that can enforce these properties for any RTO algorithm. This second part examines the practical side of the sufficient conditions for feasibility and optimality (SCFO) proposed in the first and focuses on how they may be enforced in real application, where much of the knowledge required for the conceptual SCFO is unavailable. Methods for improving convergence speed are also considered.

\;


\end{abstract}

\end{frontmatter}

\section{Review: Sufficient Conditions for Feasibility and Optimality}

The standard RTO problem requires minimizing an unknown cost function $\phi_p : \mathbb{R}^{n_u} \rightarrow \mathbb{R}$ subject to unknown inequality constraints over a relevant input space:

\begin{equation}\label{eq:realopt}
\begin{array}{l}
\mathop {{\rm{minimize}}}\limits_{\bf{u}}\hspace{4mm}\phi_p ({\bf{u}}) \\
{\rm{subject}}\hspace{1mm}{\rm{to}}\hspace{3mm}{\bf{G}}_p({\bf{u}}) \preceq {\bf{0}} \\
\hspace{18.3mm}{\bf u}^L \preceq {\bf u} \preceq {\bf u}^U
\end{array}\hspace{3mm},
\end{equation}

\noindent where ${\bf{G}}_p$ denotes the set of $n_g$ inequality constraint functions $g_{p} : \mathbb{R}^{n_u} \rightarrow \mathbb{R}$,  ${\bf{u}} \in \mathbb{R}^{n_u}$ the set of decision variables (inputs), and ${\bf u}^L$ and ${\bf u}^U$ the lower and upper limits, respectively, on these variables (i.e. the ``box'' constraints).

In the most general sense, an RTO \emph{algorithm} is some law, denoted here by $\Gamma(\cdot)$, that uses the previous and current input-output data with respect to the current iteration $k$ to calculate the optimal input target for the next iteration $k+1$:

\begin{equation}\label{eq:RTOalgo}
{\bf u}_{k+1}^* = \Gamma({\bf u}_0,...,{\bf u}_k,{\bf y}_0,...,{\bf y}_k),
\end{equation}

\noindent with ${\bf y} \in \mathbb{R}^{n_y}$ the set of measured outputs. In many cases, an input filter gain, $K_k \in [0,1]$, may be used to dampen the step taken by the RTO adaptation \citep{Brdys:05,Bunin:11b}:

\begin{equation}\label{eq:inputfilter}
{\bf{u}}_{k+1}  = {\bf{u}}_k  + K_k \left( {\bf{u}}^*_{k + 1} - {\bf{u}}_k \right),
\end{equation} 

\noindent where we have used the subscript $k$ to emphasize the fact that the gain itself may be adapted at every iteration.

As discussed in the companion paper \citep{Bunin:12b}, there are a number of requirements that any good RTO algorithm may be expected to satisfy. Primarily, the iterates generated by the algorithm should not violate any hard constraints, as doing so may cause significant equipment damage and/or promote hazardous or undesirable conditions when applied to the actual plant\footnote{Unless stated otherwise (i.e. in Section 5), we will treat all of the constraints in ${\bf G}_p$ as hard constraints.}. The iterates should also converge to a true solution of (\ref{eq:realopt}) and not to a suboptimal point. Finally, it is preferred that the algorithm yield monotonic improvement with respect to the cost, so that anyone making the choice to apply RTO may at least have the assurance that doing so will not cause their process to become \emph{less} optimal, even if only temporarily, when RTO is used.

It was shown in the previous work that such requirements could be met for the general RTO algorithm $\Gamma (\cdot)$, under reasonable assumptions, if the following sufficient conditions for feasibility and optimality (SCFO) could be made to hold at every iteration $k$ with the input filtering law (\ref{eq:inputfilter}) in place:

\begin{equation}\label{eq:suffcond2}
\nabla g_{p,j}({\bf{u}}_k)^T ({\bf{u}}_{k+1}^* - {\bf{u}}_{k}) < 0, \hspace{2mm} \forall j: g_{p,j}({\bf{u}}_k) \geq -\epsilon_j
\end{equation}

\begin{equation}\label{eq:suffcond1}
\nabla \phi_{p}({\bf{u}}_k)^T ({\bf{u}}_{k+1}^* - {\bf{u}}_{k}) < 0
\end{equation}

\begin{equation}\label{eq:suffcond3}
K_k \leq \mathop{\min}\limits_{j = 1,...,n_g} \left[ \frac{\displaystyle - g_{p,j} ({\bf u}_{k} ) }{\displaystyle \sum\limits_{i = 1}^{n_u} {\kappa_{ji} | u^*_{k+1,i} - u_{k,i} | }} \right]
\end{equation}

\begin{equation}\label{eq:suffcond4}
K_{k} < -2 \displaystyle \frac{\nabla \phi_p({\bf{u}}_k)^T ({\bf{u}}^*_{k+1} - {\bf{u}}_{k})}{({\bf{u}}^*_{k+1} - {\bf{u}}_{k})^T \overline {\bf{Q}}_\phi ({\bf{u}}^*_{k+1} - {\bf{u}}_{k})}
\end{equation}

\begin{equation}\label{eq:suffcond5}
{\bf{u}}^L \preceq {\bf{u}}_{k+1}^* \preceq {\bf{u}}^U,
\end{equation}

\noindent where $\epsilon_j > 0$ and where $\kappa_{ji}$ and $\overline {\bf Q }_\phi$ are the Lipschitz constants for the constraints and the quadratic upper bound for the cost, respectively, and must ensure the following:

\begin{equation}\label{eq:linbound}
g_{p,j} ({\bf{u}}_{k+1}) - g_{p,j} ({\bf{u}}_k) < \displaystyle \sum\limits_{i = 1}^{n_u} {\kappa_{ji} | u_{k+1,i} - u_{k,i} | }
\end{equation}

\begin{equation}\label{eq:Qbound}
\begin{array}{l}
\phi_p ({\bf{u}}_{k+1}) - \phi_p ({\bf{u}}_{k}) \leq \vspace{1mm} \\
\hspace{10mm}\nabla \phi_p ({\bf{u}}_k)^T ({\bf{u}}_{k+1} - {\bf{u}}_{k}) + \displaystyle \frac{1}{2}({\bf{u}}_{k+1} - {\bf{u}}_{k})^T \overline {\bf{Q}}_\phi ({\bf{u}}_{k+1} - {\bf{u}}_{k})
\end{array},
\end{equation}

\noindent for any ${\bf{u}}_{k}, {\bf{u}}_{k+1} \in \mathcal{I}$ ($\mathcal{I}$ being the relevant input space defined by ${\bf u}^L \preceq {\bf u} \preceq {\bf u}^U$), with ${\bf{u}}_{k} \neq {\bf{u}}_{k+1}$ for (\ref{eq:linbound}). The SCFO may be summarized qualitatively as follows:

\begin{itemize}
\item Conditions (\ref{eq:suffcond2}) and (\ref{eq:suffcond1}) ensure that the target optimum provided by the RTO algorithm, and therefore the filtered input, represents a local feasible cost descent direction. Condition (\ref{eq:suffcond2}) is only applied for constraints that are approaching activity (they are ``$\epsilon$-active'') and ensures that these constraints do not become active \emph{prematurely}, which, by Condition (\ref{eq:suffcond3}), would cause the algorithm to stop with $K_k = 0$. Not surprisingly, these conditions mimic the geometric conditions for a Karush-Kuhn-Tucker (KKT) point, in that they must be satisfied for a small enough ${\boldsymbol \epsilon}$ unless the current iterate has no local feasible cost descent direction (i.e. meets the KKT stationarity condition).
\item Condition (\ref{eq:suffcond3}) is the feasibility condition, which limits the step size of the adaptation so as to remain inside the feasible space.
\item Condition (\ref{eq:suffcond4}) is the condition on the monotonic cost decrease, which limits the step size of the adaptation so as not to ``lose'' the local descent.
\item Condition (\ref{eq:suffcond5}) simply states that the target optimum must satisfy the box constraints, and is generally assumed to be satisfied by construction (Assumption A4 in the companion work).
\end{itemize}

As no RTO scheme enforces all of the SCFO innately, an implementation that would force any algorithm to do so consists of first projecting the target given by the RTO algorithm to meet Conditions (\ref{eq:suffcond2}), (\ref{eq:suffcond1}), and (\ref{eq:suffcond5}):

\begin{equation}\label{eq:feasproj}
\begin{array}{l}
\bar {\bf{u}}^*_{k+1} = {\rm{arg}} \mathop {{\rm{minimize}}}\limits_{{\bf{u}}}\hspace{4mm}\left\| {\bf{u}} - {\bf{u}}^*_{k+1} \right\|_2^2 \vspace{1mm}  \\
\hspace{17mm}{\rm{subject}}\hspace{1mm}{\rm{to}}\hspace{3mm}\nabla g_{p,j}({\bf{u}}_k)^T ({\bf{u}} - {\bf{u}}_{k}) \leq -\delta_{g,j} \vspace{1mm} \\
\hspace{36mm}\forall j: g_{p,j}({\bf{u}}_k) \geq -\epsilon_j \vspace{1mm}\\
\hspace{36mm}\nabla \phi_{p}({\bf{u}}_k)^T ({\bf{u}} - {\bf{u}}_{k}) \leq -\delta_\phi \vspace{1mm} \\
\hspace{36mm}{\bf{u}}^L \preceq {\bf{u}} \preceq {\bf{u}}^U
\end{array},
\end{equation}

\noindent where ${\boldsymbol \delta}_{g}, \delta_{\phi} \succ {\bf 0}$ serve to approximate the strict local descent spaces\footnote{This is the ``implementable'' version of Conditions (\ref{eq:suffcond2}) and (\ref{eq:suffcond1}), and will be used by default from now on.}, and then filtering the projected target with respect to Conditions (\ref{eq:suffcond3}) and (\ref{eq:suffcond4}):

\begin{equation}\label{eq:filtertarget}
\begin{array}{l}
K_{k} := \mathop {\min} \Bigg \{ \mathop{\min}\limits_{j = 1,...,n_g} \left[\displaystyle \frac{{-g_{p,j}({\bf{u}}_{k } )}}{\displaystyle \sum\limits_{i = 1}^{n_u} {\kappa_{ji} | \bar u^*_{k+1,i} - u_{k,i} | }} \right],  \\
\hspace{25mm} -1.99 \displaystyle \frac{\nabla \phi_p({\bf{u}}_k)^T (\bar {\bf{u}}^*_{k+1} - {\bf{u}}_{k})}{(\bar {\bf{u}}^*_{k+1} - {\bf{u}}_{k})^T \overline {\bf{Q}}_\phi (\bar {\bf{u}}^*_{k+1} - {\bf{u}}_{k})} \Bigg \} \vspace{2mm}\\
K_k > 1 \rightarrow K_k := 1 \vspace{1mm} \\
{\bf{u}}_{k+1}  = {\bf{u}}_k  + K_k \left( \bar {\bf{u}}^*_{k + 1} - {\bf{u}}_k \right)
\end{array},
\end{equation}

\noindent where the heuristic choice of taking the largest step possible is made, as this generally leads to convergence in less iterations. For a detailed theoretical treatment of the SCFO, as well as of their basic implementation strategy, the interested reader is referred to the companion paper \citep{Bunin:12b} as a prerequisite for the material discussed in this one.

While the above offer very strong \emph{conceptual} guidelines as to how an RTO algorithm may be run to achieve feasible-side convergence to a KKT point, there remain a number of important implementation issues that must be addressed in order for such conditions to be applicable in practice. Specifically, the conditions as stated rely on the \emph{exact} plant gradients for both $\phi_p$ and $g_{p,j}$, as well as on knowing the appropriate values of $\kappa_{ji}$ and $\overline {\bf Q}_\phi$ so that (\ref{eq:linbound}) and (\ref{eq:Qbound}) are satisfied globally on $\mathcal{I}$. Both Conditions (\ref{eq:suffcond2}) and (\ref{eq:suffcond3}) rely on being able to measure $g_{p,j}({\bf u}_k)$ exactly. Unfortunately, many or all of these requirements are unrealistic in the practical setting.

It is the presence of these issues that motivates this second paper, whose purpose is to provide a rigorous treatment of how the SCFO may be enforced anyway despite this lack of perfect knowledge. We start by considering the subject of gradient uncertainty in Section 2, and derive a robust version of Conditions (\ref{eq:suffcond2}), (\ref{eq:suffcond1}), and (\ref{eq:suffcond4}) for the cases where \emph{bounded} gradient estimates are available. The case of noisy or estimated constraint values is then discussed in Section 3, and robust versions of Conditions (\ref{eq:suffcond2}) and (\ref{eq:suffcond3}) are proposed. Section 4 addresses the subject of upper bounds, and proposes means by which they may be chosen or refined. Section 5 then looks at the optimistic side of things, and focuses on the cases where one could make the innately conservative SCFO less so by either allowing temporary constraint violations or assuming that some parts of the RTO problem are certain and known analytically, thereby improving convergence speed. Finally, Section 6 concludes the paper with a reflection on what has been addressed and an outlook on future work.

We warn the reader that some of the results proposed are given without proofs/derivations, as the essential work for these has already been carried out in the companion article, and repeating a very similar derivation to obtain a very similar result has been deemed as unnecessary by the authors. As such, we \emph{urge the reader} to have read the previous paper prior to starting this one. All ``truly new'' results are, however, proven with the same rigor and detail as before.

\section{Gradient Uncertainty}

In the majority of applications, the concept of a \emph{gradient measurement} is nonexistent, and one must usually estimate a function's derivatives from its 0$^{th}$-order measurements (via, for example, finite differences). This is a general engineering problem with applications in multiple fields \citep{Meyer:01,Brekelmans:05,Yeow:10,Correa:11} including RTO \citep{Brdys:05,Marchetti:10,Rodger:10,Bunin:12}. It is also a challenging problem, since the almost ubiquitous presence of measurement noise, the decoupling that must be made to calculate the different directional components of the gradient, and the innate locality of the gradient all make it close to impossible to estimate a gradient exactly. As such, one must, in the best case, work with inaccurate gradient estimates and not the plant gradient itself. While rarely, if at all, treated in the RTO literature, we do note that unconstrained optimization with uncertain gradients has been looked at to a fair extent in the numerical optimization context \citep{Carter:91,Carter:93,Gratton:11}, with the general conclusion that an algorithm can be subject to moderate \emph{multiplicative} gradient noise and still be able to converge to the optimum, while the presence of additive noise cannot lead to full convergence due to resolution limits \citep{Kelley:03}.

In this section, we consider the problem of satisfying the SCFO when the estimated plant gradients are inaccurate, and show that this may still be possible when bounds on the (additive) gradient estimation error are available and serve to define an uncertainty set to which the plant gradient must belong. Furthermore, it is demonstrated how knowledge of such bounds may be used to generate robust versions of both Projection (\ref{eq:feasproj}) and the filter gain limit (\ref{eq:suffcond4}). An implementation algorithm to enforce SCFO satisfaction in the presence of gradient uncertainty is proposed, and the results of a large numerical study are given.

\subsection{SCFO with Bounded Gradient Estimates}

As the most basic approach, which we will call the \emph{0-robustness} projection, we could attempt to use the estimated gradients, denoted by $ \nabla \hat g_{p,j}({\bf{u}}_k)$ and $\nabla \hat \phi_{p}({\bf{u}}_k)$, directly in Projection (\ref{eq:feasproj}). However, such an approach would be risky at best, since the estimated gradients are likely to be erroneous and what would be satisfied would no longer be the true SCFO but some corrupted version that may not give the desired guarantees. The issues with using gradient estimates directly may be summed up as follows:

\begin{itemize}
\item Inability to stay away from $\epsilon$-active constraints, since the local descent space generated by the gradient estimate may not be correct, and $\nabla \hat g_{p,j}({\bf{u}}_k)^T$ $ ({\bf{u}} - {\bf{u}}_{k}) \leq -\delta_{g,j} \not \Rightarrow \nabla g_{p,j}({\bf{u}}_k)^T ({\bf{u}} - {\bf{u}}_{k}) \leq -\delta_{g,j}$.
\item Inability to decrease the cost at every iteration due to errors both in (\ref{eq:suffcond1}) and (\ref{eq:suffcond4}).
\item Convergence to a false KKT point when the projection becomes infeasible for a false set of gradients, premature convergence, or no convergence at all.
\end{itemize}

It may, however, be possible to satisfy the true SCFO when the true plant gradients may be identified as belonging to a compact uncertainty set:

\begin{equation}\label{eq:boundedgrad}
\begin{array}{l}
\nabla \underline g_{p,j}({\bf{u}}_k) \preceq \nabla g_{p,j}({\bf{u}}_k) \preceq \nabla \overline g_{p,j}({\bf{u}}_k), \hspace{2mm}\forall j: g_{p,j}({\bf{u}}_k) \geq -\epsilon_j \\
\nabla \underline \phi_{p}({\bf{u}}_k) \preceq \nabla \phi_{p}({\bf{u}}_k) \preceq \nabla \overline \phi_{p}({\bf{u}}_k)
\end{array}.
\end{equation}

\noindent While not necessarily easy to obtain, bounds on the gradients are usually possible when given an additional assumption on the corresponding function. In finite-difference schemes, for example, one may bound the estimation error due to noise and to truncation by considering the worst-case amplitude of the noise together with the worst-case evolution (i.e. the quadratic upper bound) of the function of interest \citep{Brekelmans:05,Marchetti:10}. One could also estimate bounds by assuming a certain local structure for the function \citep{Bunin2013a} -- for example, the assumption of local quasiconvexity or quasiconcavity constrains the gradients to lie in a polyhedron \citep{Bunin:12}. When such assumptions are either impossible or undesired, a less rigorous but nevertheless sensible ``general-purpose'' approach may be to take the estimated gradient and to generate an uncertainty region around it,

\begin{equation}\label{eq:genbound}
\begin{array}{l}
\nabla \overline g_{p,j}({\bf{u}}_k), \nabla \underline g_{p,j}({\bf{u}}_k) = \nabla \hat g_{p,j} ({\bf u}_k) \pm m \sigma_{\nabla,j}, \hspace{2mm}\forall j: g_{p,j}({\bf{u}}_k) \geq -\epsilon_j \\
\nabla \overline \phi_{p}({\bf{u}}_k), \nabla \underline \phi_{p}({\bf{u}}_k) = \nabla \hat \phi_{p} ({\bf u}_k) \pm m \sigma_{\nabla,\phi}
\end{array},
\end{equation}

\noindent by increasing the value of the positive scalar $m$, with $\sigma_\nabla \in \mathbb{R}^{n_u}_+$ denoting different scaling vectors. Finally, note that for the constraints these bounds simply reduce, in the worst case, to the Lipschitz constants:  

\begin{equation}\label{eq:lip0}
-\kappa_{ji} < \frac{\partial g_{p,j}}{\partial u_i} \Big |_{\bf{u}} < \kappa_{ji}, \;\; \forall {\bf u} \in \mathcal{I}.
\end{equation}

In the discussion that follows, we do not question the method used to obtain the bounds, but assume that the real gradient lies, almost surely, within them. With the availability of this additional knowledge, the \emph{robust} versions of (\ref{eq:suffcond2}) and (\ref{eq:suffcond1}) may then be written as:

\begin{equation}\label{eq:suffconrelax}
\begin{array}{l}
\nabla \tilde g_{p,j}({\bf{u}}_k)^T ({\bf{u}}^*_{k+1} - {\bf{u}}_{k}) \leq -\delta_{g,j} \vspace{1mm}\\ 
\hspace{15mm}\forall \nabla \tilde g_{p,j}({\bf{u}}_k) \in [\nabla \underline g_{p,j}({\bf{u}}_k), \nabla \overline g_{p,j}({\bf{u}}_k)], \hspace{2mm}\forall j: g_{p,j}({\bf{u}}_k) \geq -\epsilon_j \vspace{1mm} \\
\nabla \tilde \phi_{p}({\bf{u}}_k)^T ({\bf{u}}^*_{k+1} - {\bf{u}}_{k}) \leq -\delta_\phi \vspace{1mm}\\
\hspace{15mm}\forall \nabla \tilde \phi_{p}({\bf{u}}_k) \in [\nabla \underline \phi_{p}({\bf{u}}_k), \nabla \overline \phi_{p}({\bf{u}}_k)] 
\end{array},
\end{equation}

\noindent where it is evident that satisfying the above will naturally satisfy the true SCFO, since the true SCFO constraint inequalities are included in the set generated by (\ref{eq:suffconrelax}). This is the typical strategy in robust analysis, as we simply try to satisfy the SCFO for all of the possible scenarios in the uncertain set.

As stated, (\ref{eq:suffconrelax}) leads to a semi-infinite set of constraints and is generally not tractable numerically if incorporated into a projection like (\ref{eq:feasproj}). We now propose an equivalent, tractable version for a given constraint $g_{p,j}$, noting that the exact same analysis applies to the other cost and constraint gradients.

\begin{lem}{\bf (Semi-Infinite to Finite Reformulation)}

Consider the semi-infinite set of constraints:

\begin{equation}\label{eq:seminf}
\begin{array}{l}
\nabla \tilde g_{p,j}({\bf{u}}_k)^T ({\bf{u}}^*_{k+1} - {\bf{u}}_{k}) \leq -\delta_{g,j}, \; \forall \nabla \tilde g_{p,j}({\bf{u}}_k) \in [\nabla \underline g_{p,j}({\bf{u}}_k), \nabla \overline g_{p,j}({\bf{u}}_k)]
\end{array},
\end{equation}

\noindent and the reformulation:

\begin{equation}\label{eq:reform}
\begin{array}{l}
\displaystyle \sum\limits_{i = 1}^{n_u} {s_{ji}} \leq -\delta_{g,j} \vspace{1mm} \\
\displaystyle \frac{\partial \underline g_{p,j}}{\partial u_i} \Big | _{{\bf{u}}_k} (u_{k+1,i}^* - u_{k,i} ) \leq s_{ji} \vspace{1mm}\\
\displaystyle \frac{\partial \overline g_{p,j}}{\partial u_i} \Big | _{{\bf{u}}_k} (u_{k+1,i}^* - u_{k,i} ) \leq s_{ji} \vspace{1mm}\\
\end{array},
\end{equation}

\noindent where $s_{ji}$ denote auxiliary slack variables. (\ref{eq:seminf}) and (\ref{eq:reform}) are equivalent in the following sense:

\begin{equation}\label{eq:equiv}
(\ref{eq:seminf}) \; {\rm feasible} \; {\rm for} \; {\bf u}_{k+1}^* \Leftrightarrow (\ref{eq:reform}) \; {\rm feasible} \; {\rm for} \; {\bf u}_{k+1}^* \; {\rm and} \; {\rm some} \; s_{ji}.
\end{equation}

\end{lem}

\begin{pf}

For convenience, we rewrite (\ref{eq:seminf}) in its summation form:

\begin{equation}\label{eq:der1}
\displaystyle \sum\limits_{i = 1}^{n_u} {\frac{\partial \tilde g_{p,j}}{\partial u_i} \Big | _{{\bf{u}}_k} (u^*_{k+1,i} - u_{k,i} ) } \leq -\delta_{g,j}.
\end{equation}

Owing to the box nature of the gradient uncertainty set, it is easy to see that the worst-case value of each individual component may be bounded as follows:

\begin{equation}\label{eq:der2}
\begin{array}{l}
u^*_{k+1,i} - u_{k,i} \geq 0 \Rightarrow \vspace{2mm} \\
\displaystyle {\frac{\partial \tilde g_{p,j}}{\partial u_i} \Big | _{{\bf{u}}_k} (u^*_{k+1,i} - u_{k,i} ) } \leq {\frac{\partial \overline g_{p,j}}{\partial u_i} \Big | _{{\bf{u}}_k} (u^*_{k+1,i} - u_{k,i} ) } \vspace{2mm} \\
u^*_{k+1,i} - u_{k,i} \leq 0 \Rightarrow \vspace{2mm} \\
\displaystyle {\frac{\partial \tilde g_{p,j}}{\partial u_i} \Big | _{{\bf{u}}_k} (u^*_{k+1,i} - u_{k,i} ) } \leq {\frac{\partial \underline g_{p,j}}{\partial u_i} \Big | _{{\bf{u}}_k} (u^*_{k+1,i} - u_{k,i} ) }
\end{array}.
\end{equation}

We first suppose that there exist ${\bf u}_{k+1}^*, s_{ji}$ so that (\ref{eq:reform}) is satisfied. From (\ref{eq:der2}) and the last two constraints of (\ref{eq:reform}), it follows that:

\begin{equation}\label{eq:th1D}
\begin{array}{l}
\displaystyle \frac{\partial \tilde g_{p,j}}{\partial u_i} \Big | _{{\bf{u}}_k} (u_{k+1,i}^* - u_{k,i} ) \leq \\
\hspace{10mm} \displaystyle {{\rm{max}} \left({\frac{\partial \overline g_{p,j}}{\partial u_i} \Big | _{{\bf{u}}_k} (u_{k+1,i}^* - u_{k,i} ) }, {\frac{\partial \underline g_{p,j}}{\partial u_i} \Big | _{{\bf{u}}_k} (u_{k+1,i}^* - u_{k,i} ) }\right)} \leq s_{ji}
\end{array},
\end{equation}

\noindent which, given the first constraint of (\ref{eq:reform}), makes it clear that:

\begin{equation}\label{eq:th1E}
\displaystyle \sum\limits_{i = 1}^{n_u} \frac{\partial \tilde g_{p,j}}{\partial u_i} \Big | _{{\bf{u}}_k} (u_{k+1,i}^* - u_{k,i} ) = \nabla \tilde g_{p,j} ({\bf u}_k)^T ({\bf u}_{k+1}^* - {\bf u}_k) \leq -\delta_{g,j}.
\end{equation}

We have thereby proven that:

\begin{equation}\label{eq:equivA}
(\ref{eq:seminf}) \; {\rm feasible} \; {\rm for} \; {\bf u}_{k+1}^* \Leftarrow (\ref{eq:reform}) \; {\rm feasible} \; {\rm for} \; {\bf u}_{k+1}^* \; {\rm and} \; {\rm some} \; s_{ji}.
\end{equation}

Now, suppose the converse where there exists some ${\bf u}_{k+1}^*$ so that (\ref{eq:seminf}) is satisfied. Clearly, this must also hold for the cases where:

\begin{equation}\label{eq:th1H}
\displaystyle \sum\limits_{i = 1}^{n_u} {{\rm{max}} \left({\frac{\partial \overline g_{p,j}}{\partial u_i} \Big | _{{\bf{u}}_k} (u_{k+1,i}^* - u_{k,i} ) }, {\frac{\partial \underline g_{p,j}}{\partial u_i} \Big | _{{\bf{u}}_k} (u_{k+1,i}^* - u_{k,i} ) }\right)} \leq -\delta_{g,j},
\end{equation}

\noindent since all of the elements in the summation and in the maximum operator are in the gradient uncertainty set.

Choosing

\begin{equation}\label{eq:th1I}
s_{ji} =  {{\rm{max}} \left({\frac{\partial \overline g_{p,j}}{\partial u_i} \Big | _{{\bf{u}}_k} (u_{k+1,i}^* - u_{k,i} ) }, {\frac{\partial \underline g_{p,j}}{\partial u_i} \Big | _{{\bf{u}}_k} (u_{k+1,i}^* - u_{k,i} ) }\right)}
\end{equation}

\noindent clearly satisfies (\ref{eq:reform}), thereby proving that:

\begin{equation}\label{eq:equivB}
(\ref{eq:seminf}) \; {\rm feasible} \; {\rm for} \; {\bf u}_{k+1}^* \Rightarrow (\ref{eq:reform}) \; {\rm feasible} \; {\rm for} \; {\bf u}_{k+1}^* \; {\rm and} \; {\rm some} \; s_{ji},
\end{equation}

\noindent and, consequently, (\ref{eq:equiv}). \qed

\end{pf}

\begin{figure}
\begin{center}
\includegraphics[width=12cm]{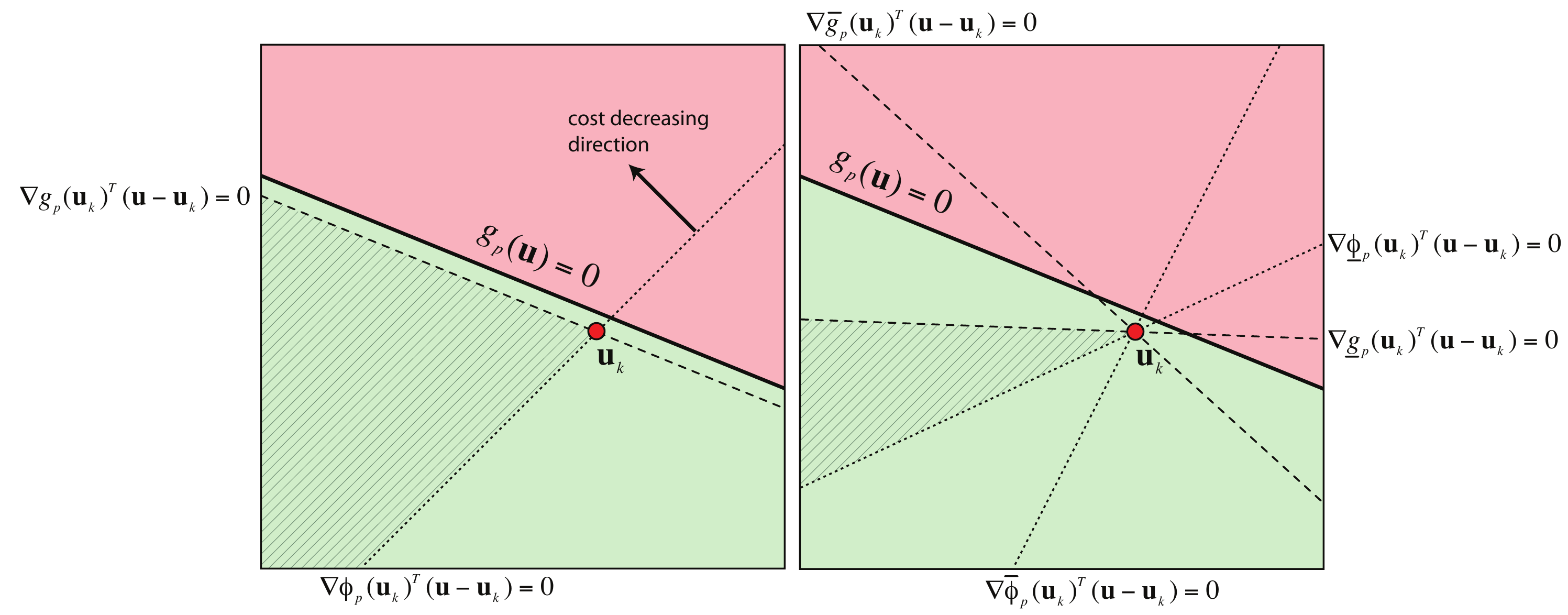}    
\caption{The use of bounded gradients reduces the feasible space for the projection problem. The nominal problem with certain gradients is given on the left, where a local descent space for both the cost and constraint is sought and corresponds to the intersection of two descent halfspaces. On the right, bounded gradients are used, which has the effect of compressing the halfspaces into cones and reducing the overall feasible space of their intersection. Lined regions denote the feasible areas in both cases, with green and red regions denoting the feasible and infeasible spaces of the RTO problem, respectively.}
\label{fig:graduncertain}
\end{center}
\end{figure}

Replacing the local descent constraints in Projection (\ref{eq:feasproj}) with the reformulation of (\ref{eq:reform}) results in the robust projection:

\begin{equation}\label{eq:feasprojrob}
\begin{array}{l}
\bar {\bf{u}}^*_{k+1} = {\rm{arg}}\mathop {{\rm{minimize}}}\limits_{{\bf{u}},{\bf{S}},{\bf{s}}_\phi}\hspace{4mm}\left\| {\bf{u}} - {\bf{u}}^*_{k+1} \right\|_2^2 \vspace{1mm}  \\
\hspace{17mm}{\rm{subject}}\hspace{1mm}{\rm{to}}\hspace{3mm}\displaystyle \sum\limits_{i = 1}^{n_u} {s_{ji}} \leq -\delta_{g,j} \vspace{1mm} \\
\hspace{36mm}\displaystyle \frac{\partial \underline g_{p,j}}{\partial u_i} \Big | _{{\bf{u}}_k} (u_{i} - u_{k,i} ) \leq s_{ji} \vspace{1mm}\\
\hspace{36mm}\displaystyle \frac{\partial \overline g_{p,j}}{\partial u_i} \Big | _{{\bf{u}}_k} (u_{i} - u_{k,i} ) \leq s_{ji} \vspace{1mm}\\
\hspace{36mm}\forall j: g_{p,j}({\bf{u}}_k) \geq -\epsilon_j \vspace{1mm}\\
\hspace{36mm}\displaystyle \sum\limits_{i = 1}^{n_u} {s_{\phi,i}} \leq -\delta_\phi \vspace{1mm} \\
\hspace{36mm}\displaystyle \frac{\partial \underline \phi_{p}}{\partial u_i} \Big | _{{\bf{u}}_k} (u_{i} - u_{k,i} ) \leq s_{\phi,i} \vspace{1mm}\\
\hspace{36mm}\displaystyle \frac{\partial \overline \phi_{p}}{\partial u_i} \Big | _{{\bf{u}}_k} (u_{i} - u_{k,i} ) \leq s_{\phi,i} \vspace{1mm}\\
\hspace{36mm}{\bf{u}}^L \preceq {\bf{u}} \preceq {\bf{u}}^U
\end{array},
\end{equation}

\noindent with ${\bf{S}} \in \mathbb{R}^{n_g \times n_u}$ the matrix of the slack variables for the constraints with entries $s_{ji}$, and ${\bf{s}}_\phi \in \mathbb{R}^{n_u}$ the vector of slack variables, $s_{\phi,i}$, corresponding to the cost. Though more involved than (\ref{eq:feasproj}), this is still a quadratic program with, at most, $n_u (n_g + 2)$ variables and $2n_u (n_g + 2)$ linear constraints. We provide a geometrical interpretation of incorporating gradient uncertainty and the effect that it has on the feasible space of the projection (i.e. its reduction) in Figure \ref{fig:graduncertain}.

We now analyze Condition (\ref{eq:suffcond4}), and write its robust version as:

\begin{equation}\label{eq:filtercostrob}
\begin{array}{l}
K_{k} < -2 \displaystyle \frac{\nabla \tilde \phi_p({\bf{u}}_k)^T ( \bar {\bf{u}}^*_{k+1} - {\bf{u}}_{k})}{( \bar {\bf{u}}^*_{k+1} - {\bf{u}}_{k})^T \overline {\bf{Q}}_\phi (\bar {\bf{u}}^*_{k+1} - {\bf{u}}_{k})} \vspace{2mm}\\
\hspace{20mm}\forall \nabla \tilde \phi_{p}({\bf{u}}_k) \in [\nabla \underline \phi_{p}({\bf{u}}_k), \nabla \overline \phi_{p}({\bf{u}}_k)]
\end{array}.
\end{equation}

As $\nabla \tilde \phi_p({\bf{u}}_k)^T ( \bar {\bf{u}}^*_{k+1} - {\bf{u}}_{k})$ is forced to be negative, we are only interested in its maximum value (as this would result in the smallest upper bound for $K_k$). From (\ref{eq:th1D}), we may define the robust version of this condition as:

\begin{equation}\label{eq:costfilter2}
K_k < -2\frac{\displaystyle \sum\limits_{i = 1}^{n_u} {{\rm{max}} \left({\frac{\partial \overline \phi_p}{\partial u_i} \Big | _{{\bf{u}}_k} (\bar u^*_{k+1,i} - u_{k,i} ) }, {\frac{\partial \underline \phi_p}{\partial u_i} \Big | _{{\bf{u}}_k} (\bar u^*_{k+1,i} - u_{k,i} ) }\right)}}{(\bar {\bf{u}}^*_{k+1} - {\bf{u}}_{k})^T \overline {\bf{Q}}_\phi (\bar {\bf{u}}^*_{k+1} - {\bf{u}}_{k})}.
\end{equation}

\subsection{Implementation, Infeasibility, and Partial Robustness}

By Theorem 5 in the companion paper, the robust projection is conceptually ``clean'' in that it preserves all of the desired convergence guarantees by satisfying the true plant SCFO despite not having accurate gradient values. Practically, however, it is not difficult to show that there are many cases for which enforcing the SCFO for the entire gradient uncertainty set is simply impossible (returning to Figure \ref{fig:graduncertain}, one may imagine the feasible region disappearing entirely as the uncertainty in the gradients is increased). This rather natural shortcoming is formalized in the following theorem.

\begin{thm}{\bf (Infeasibility of the Robust Projection)}

Let $\mathcal{M}_k$ be a matrix with a semi-infinite number of rows and $n_u$ columns representing the entire semi-infinite set of constraints in (\ref{eq:suffconrelax}), along with any box constraints that are active, so that any ${\bf u}_{k+1}^*$ satisfying

\begin{equation}\label{eq:M}
\begin{array}{l}
\mathcal{M}_k ({\bf{u}}_{k+1}^* - {\bf{u}}_{k}) \preceq \left[ {\begin{array}{*{20}c}
   -\boldsymbol{\delta}_{g,\mathcal{M}}  \\
   -\boldsymbol{\delta}_{\phi,\mathcal{M}} \\
   {\bf{0}} \\
\end{array}} \right]
\end{array}
\end{equation}

\noindent satisfies the SCFO robustly. For a fixed set of ${\boldsymbol \epsilon}$, it follows that Projection (\ref{eq:feasprojrob}) is infeasible for all choices of $\boldsymbol{\delta}_{g,\mathcal{M}}, \boldsymbol{\delta}_{\phi,\mathcal{M}} \succ {\bf 0}$ iff there exist one or more rows of $\mathcal{M}_k$ whose scaled sum is equal to ${\bf 0}$ for some choice of strictly positive scalars.

\end{thm}

\begin{pf}

A similar analysis has already been carried out in the companion work for the case where $\mathcal{M}_k$ was finite-dimensional (see Theorems 1, 4, and 5 therein). Ready analogues to Theorem 1 (of the companion work) exist for the semi-infinite case \citep{Goberna:84} and state that the system (\ref{eq:M}) is infeasible iff a negative spanning exists among some of its members (as in the finite case). By the equivalence proven in Lemma 1, we know that an infeasibility in (\ref{eq:M}) implies an infeasibility in the corresponding finite reformulation, and thereby in Projection (\ref{eq:feasprojrob}). \qed

\end{pf}

The above has an interesting practical interpretation. First, note that, like in the ideal case (see Theorem 5 in the companion paper), there are some cases for which the projection is certainly infeasible:

\begin{itemize}
\item $\nabla \phi_p ({\bf u}_k) = {\bf 0}$,
\item $\exists j : \nabla g_{p,j} ({\bf u}_k) = {\bf 0}, \;\; g_{p,j} ({\bf u}_k) \geq -\epsilon_j $,
\item $\exists {\boldsymbol \mu} \in \mathbb{R}^{n_g}_+ : \nabla \phi_p ({\bf u}_k) + \displaystyle \sum\limits_{j = 1}^{n_g} \mu_j \nabla g_{p,j} ({\bf u}_k) = {\bf 0}$,
\end{itemize}

\noindent where the first and third cases correspond to stationary KKT points and the second is a pathological case that we assume to be of little practical concern (see Assumption A3 in the companion paper). As the matrix $\mathcal{M}_k$ includes the plant gradients, it follows that the robust projection is infeasible as each of these three cases implies a negative spanning. Theorem 1 says more than this, however, as the robust projection is guaranteed to be infeasible \emph{whenever} gradients corresponding to the above cases are present in $\mathcal{M}_k$, regardless or not of whether they correspond to a plant KKT point.

There is both an optimistic and a pessimistic way to interpret this property. For the former, we can make the argument that the infeasibility of the robust projection implies, in some sense, that we are already in the neighborhood of a stationary point. To see why this is true, consider being at a point that is far from a stationary point and for which the gradient bounds are reasonable (i.e. do not deviate too far from the true gradient). It is then unlikely that any of the gradients in the uncertainty set satisfy a stationarity condition since the plant gradients do not by a significant amount. As one then approaches a stationary point and the plant gradients approach stationarity, the chance of small deviations from the plant gradients (i.e. all of the gradients in the uncertainty set) satisfying the stationarity condition increases. The pessimistic view is that the size of such a neighborhood is unknown and can become arbitrarily large as the quality of the bounds deteriorates (it is evident that we cannot, for example, ensure local descent for \emph{every} possible gradient in $\mathbb{R}^{n_u}$ when the bounds are arbitrarily loose, irrespective of ${\bf u}$).

\begin{figure}
\begin{center}
\includegraphics[width=8cm]{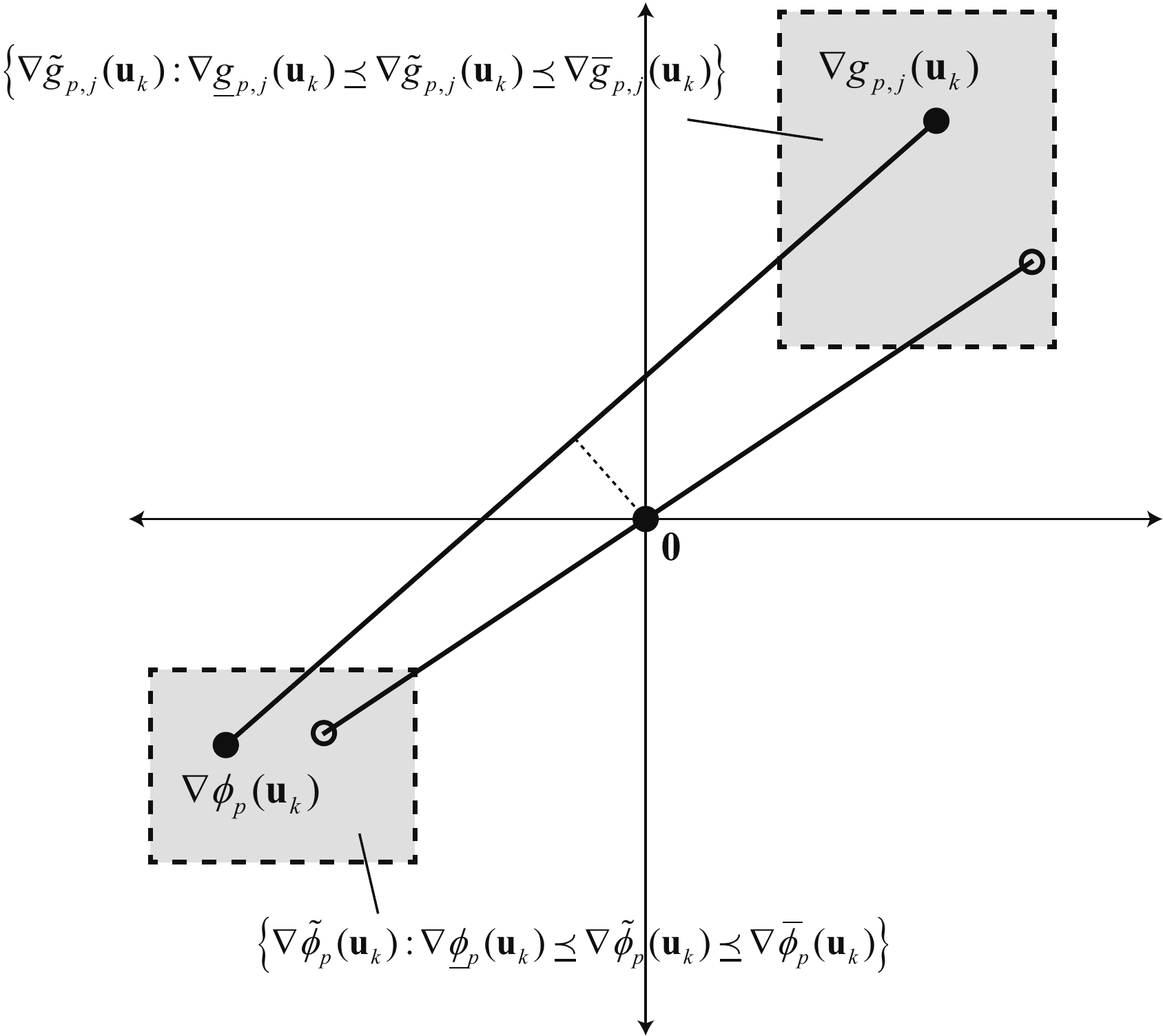}    
\caption{The relationship between the uncertainty in the gradients and the infeasibility of the robust projection due to ``false stationarity''. The axes here may be interpreted as the values of the derivatives in the different directions.}
\label{fig:graduncert}
\end{center}
\end{figure}

A geometrical illustration of these concepts is given in Figure \ref{fig:graduncert}. Here, a simple example with a single constraint is used, with the true gradients of both the cost and constraint being plotted with the filled, black circles. It is clear that there is no negative spanning between the two as their scaled sums (denoted by the line that connects them) does not pass through the origin, and the distance between this line and the origin (given by the dotted line) may be viewed as a sort of KKT metric. When the gradients are uncertain, however, and when the true gradients are known only to lie somewhere in the uncertainty set (the two boxes), it becomes easy to find a pair of estimates (given by the hollow points) that satisfy the stationarity condition and whose positively scaled sum is ${\bf 0}$. The result is that it is impossible to satisfy the robust SCFO in this case, since the robust SCFO must account for all possible gradient combinations but cannot be satisfied if a certain combination meets the stationarity condition. It should be evident that this undesired effect may only be eliminated by tightening or shrinking the uncertainty boxes.

This leaves us with the question of what should be done when the projection becomes infeasible, as it inevitably will. There are, in our opinion, three options:

\begin{itemize}
\item Stop RTO adaptation and declare convergence to a neighborhood of a stationary point.
\item Temporarily stop RTO adaptation and perform local perturbations to get a tighter estimate of the local gradients (via, e.g., a finite-difference approach). This serves to shrink the gradient uncertainty set robustly, and may be done until the projection becomes feasible. 
\item Reduce the robustness bounds, i.e. shrink the uncertainty set artificially, until the projection becomes feasible. This approach is \emph{not} robust, but sacrifices as little robustness as possible in order to make the projection feasible. As such, this may be viewed as a \emph{partial-robustness} (or stochastic, if probability distributions on the uncertainty are assumed) approach.
\end{itemize}

It should be clear that the first choice has the least attractive properties for cases where the process in question is expected to run for a very large number of iterations, as the suboptimality accepted by simply stopping is likely to add up to become significant over the entire operation time. By contrast, the second approach is well-suited for such cases since it still guarantees finite-time convergence (assuming that it takes a finite number of perturbation experiments to shrink the gradient uncertainty set sufficiently), although it does not guarantee monotonic cost improvement as these perturbations may lead to cost increases. The third approach may be appropriate when fast cost improvement is desired without rigorous satisfaction of the SCFO -- no additional experiments to shrink the uncertainty set are required and a certain degree of robustness is always retained. Note that all three may be used in the course of the operation of a single process, however -- for a problem that is initialized with a significantly suboptimal point, one could begin with the third strategy, then shift to the second once significant improvements have been attained, and then adopt the first once improvements become negligible.

We finish this discussion by proposing the following modified algorithm for enforcing the SCFO in the case of gradient uncertainty, focusing specifically on the third approach (extensions to using the first or second methods exclusively may be obtained by trivial modifications of either this algorithm or the one presented in the companion paper):

\vspace{2mm}
\noindent {\bf{Initialization -- Done Prior to RTO}}
\vspace{2mm}

\begin{enumerate}
\item Scale the constraints and cost with respect to their ranges (e.g. if $\mathop{\min}\limits_{{\bf{u}} \in \mathcal{I}} g_p({\bf{u}}) = -100$, the scaled constraint may be re-defined as $g_p({\bf{u}}) := 0.01g_p({\bf{u}})$). 

Setting upper and lower limits on the projection parameters, let $\overline {\boldsymbol \epsilon} = \overline {\boldsymbol \delta}_g = {\bf{1}}$, $\overline \delta_\phi = 1$, and choose $\underline {\boldsymbol \epsilon}$, $\boldsymbol{\underline \delta}_g$, and $\underline \delta_\phi$ to be sufficiently small (e.g. 10$^{-6}$). 
\end{enumerate}

\vspace{2mm}
\noindent {\bf{Search for a Feasible Direction -- Before Each RTO Iteration}}
\vspace{2mm}

\begin{enumerate}
\setcounter{enumi}{1}
\item Set ${\boldsymbol \epsilon} := \overline {\boldsymbol \epsilon}$, $\boldsymbol{\delta}_g := \boldsymbol{\overline \delta}_g$, and $\delta_\phi := \overline \delta_\phi$. Set $P := 1$.
\item Check the feasibility of the 0-robustness projection -- (\ref{eq:feasproj}) with estimated gradients $\nabla \hat g_{p,j} ({\bf u}_k)$ and $\nabla \hat \phi_{p} ({\bf u}_k)$ in place of the plant gradients $\nabla g_{p,j} ({\bf u}_k)$ and $\nabla \phi_{p} ({\bf u}_k)$ -- for the given choice of ${\boldsymbol \epsilon}$, $\boldsymbol{\delta}_g$, and $\delta_\phi$.

If the 0-robustness does not have a solution, set ${\boldsymbol \epsilon} := 0.5 {\boldsymbol \epsilon}$, $\boldsymbol{\delta}_g := 0.5 \boldsymbol{\delta}_g$, $\delta_\phi := 0.5 \delta_\phi$, and attempt to re-solve the problem. If ${\boldsymbol \epsilon} \prec \underline {\boldsymbol \epsilon}$, $\boldsymbol{\delta}_g \prec \boldsymbol{\underline \delta}_g$, and $\delta_\phi < \underline \delta_\phi$, but the projection is still infeasible, then either (i) declare convergence and go to Step 5 or (ii) apply a local perturbation to the system to refine the gradient estimate and the corresponding uncertainty set, and then return to Step 2. If the 0-robustness problem is feasible, proceed to Step 4.

\item Check the feasibility of the $P$-robustness projection for the resulting choice of ${\boldsymbol \epsilon}$, $\boldsymbol{\delta}_g$, and $\delta_\phi$, where $P$ defines the gradient uncertainty set in the following manner:

\begin{equation}\label{eq:boundslack}
\begin{array}{l}
\nabla \overline g_{p,j}({\bf{u}}_k) := \nabla \hat g_{p,j} ({\bf u}_k) + P \Big ( \nabla \overline g_{p,j}({\bf{u}}_k) - \nabla \hat g_{p,j} ({\bf u}_k)  \Big ) \\
\nabla \underline g_{p,j}({\bf{u}}_k) := \nabla \hat g_{p,j} ({\bf u}_k) + P \Big ( \nabla \underline g_{p,j}({\bf{u}}_k) - \nabla \hat g_{p,j} ({\bf u}_k)  \Big ) \\
\nabla \overline \phi_{p}({\bf{u}}_k) := \nabla \hat \phi_{p} ({\bf u}_k) + P \Big (  \nabla \overline \phi_{p}({\bf{u}}_k) -  \nabla \hat \phi_{p} ({\bf u}_k) \Big ) \\
\nabla \underline \phi_{p}({\bf{u}}_k) := \nabla \hat \phi_{p} ({\bf u}_k) + P \Big (  \nabla \underline \phi_{p}({\bf{u}}_k) -  \nabla \hat \phi_{p} ({\bf u}_k) \Big )
\end{array},
\end{equation}

\noindent where the upper and lower bounds on the right-hand side are \emph{always} the original robust bounds. If the $P$-robustness projection does not have a solution, set $P := P - 0.05$ and attempt to re-solve the projection. Otherwise, solve the $P$-robust projection (\ref{eq:feasprojrob}) for the bounds defined in (\ref{eq:boundslack}) to obtain $\bar {\bf{u}}^*_{k+1}$.
\end{enumerate}

\vspace{2mm}
\noindent {\bf{Termination -- Declared Convergence to KKT Point}}
\vspace{2mm}

\begin{enumerate}
\setcounter{enumi}{4}
\item If ${\boldsymbol \epsilon} \prec \underline {\boldsymbol \epsilon}$, $\boldsymbol{\delta}_g \prec \boldsymbol{\underline \delta}_g$, $\delta_\phi < \underline \delta_\phi$, terminate, with $\bar {\bf{u}}^*_{k+1} := {\bf{u}}_{k}$.
\end{enumerate}

\vspace{2mm}

Implicit in the above algorithm is the notion that the uncertainty set reduces to the gradient estimate in the 0-robustness case. This make sense from a statistical perspective, as the estimate will usually be the best guess of the gradient in some statistical sense. Finally, we note that the reduced bounds (\ref{eq:boundslack}) should be used instead of the original bounds in the filter gain bound of (\ref{eq:costfilter2}).

\subsection{Examples}

We return to the RTO problem that was used as the illustrative example in Section 4 of the companion paper:

\begin{equation}\label{eq:prob1}
\begin{array}{l}
\mathop {{\rm{minimize}}}\limits_{u_1,u_2} \hspace{3mm} \phi_{p}({\bf{u}}) = (u_1-0.5)^2 + (u_2-0.4)^2 \\
{\rm{subject}}\hspace{1mm}{\rm{to}}\hspace{3mm}g_{p,1}({\bf{u}})= -6u^2_1 - 3.5u_1 + u_2 -0.6 \le 0 \vspace{1mm} \\
\hspace{18mm}g_{p,2}({\bf{u}})= 2u^2_1 + 0.5u_1 + u_2 -0.75 \le 0 \vspace{1mm} \\
\hspace{18mm}g_{p,3}({\bf{u}})= -u^2_1 - (u_2-0.15)^2 +0.01 \le 0 \vspace{1mm} \\
\hspace{18mm}u_1 \in [-0.5, 0.5], u_2 \in [0, 0.8]
\end{array},
\end{equation}

\noindent and has both a ``difficult'' and an ``easy'' scenario, which correspond to the initial starting points of ${\bf{u}}_0 = [-0.5, 0.05]$ and ${\bf{u}}_0 = [0, 0.4]$. We will refer to these as Problem A and Problem B, respectively.

As in the companion work, we consider a number of algorithms, and refer the reader to the first paper for the details regarding the basic implementation (i.e. the choice of projection parameters, Lipschitz constants, and algorithm definitions). In this section, we will essentially repeat the trials of the previous work but will work with noisy gradient estimates\footnote{We use $\mathcal{U}[-1,1]$ directly in the equations to denote a random number drawn from $\mathcal{U}[-1,1]$.}:

\begin{equation}\label{eq:noisygrad}
\begin{array}{l}
\displaystyle \frac{\partial \hat g_{p,j}}{\partial u_i} \Big |_{{\bf u}_k} = \frac{\partial g_{p,j}}{\partial u_i}\Big |_{{\bf u}_k} + \sigma \kappa_{ji} \mathcal{U}[-1,1] \\
\displaystyle \frac{\partial \hat \phi_{p}}{\partial u_i} \Big |_{{\bf u}_k} = \frac{\partial \phi_{p}}{\partial u_i}\Big |_{{\bf u}_k} + \sigma \kappa_{\phi,i} \mathcal{U}[-1,1]
\end{array},
\end{equation}

\noindent i.e. each derivative is corrupted by additive uniform noise that is proportional, via the adjustable scalar $\sigma$, to its Lipschitz constant. For the cost, we choose the constants ${\boldsymbol \kappa}_\phi = [2.2\;\; 0.35]^T$.

Defining the performance measure of an algorithm as the sum of its optimality losses:

\begin{equation}\label{eq:perfmetric}
L = \displaystyle \sum\limits_{i = 0}^{k_f} \left[ \phi_p ({\bf u}_i) - \phi_p ({\bf u}^*) \right] ,
\end{equation}

\noindent with $k_f$ the final iteration (1,000 for A and 100 for B) and ${\bf u}^*$ the plant optimum, we consider the following four implementations:

\begin{enumerate}[I.]
\item The nominal case (0-robustness with $\sigma = 0$).
\item The 0-robustness case, where the noisy gradients are used directly in the projection with $\sigma > 0$.
\item The robust-projection case where the magnitude of $\sigma$ is assumed to be known and the following bounds are adopted:

\begin{equation}\label{eq:exbound}
\begin{array}{l}
\nabla \hat g_{p,j}({\bf{u}}_k) - \sigma {\boldsymbol \kappa}_j \preceq \nabla g_{p,j}({\bf{u}}_k) \preceq \nabla \hat g_{p,j}({\bf{u}}_k) + \sigma {\boldsymbol \kappa}_j \\
\nabla \hat \phi_{p}({\bf{u}}_k) - \sigma {\boldsymbol \kappa}_\phi \preceq \nabla \phi_{p}({\bf{u}}_k) \preceq \nabla \hat \phi_{p}({\bf{u}}_k) + \sigma {\boldsymbol \kappa}_\phi
\end{array}.
\end{equation}

Partial robustness is accepted by reducing the bounds as in (\ref{eq:boundslack}) when the robust projection is infeasible.

\item The robust-projection case where no knowledge of $\sigma$ is assumed, a general uncertainty of $\sigma_\nabla = {\bf 1}$ is chosen, and the general-purpose bounds (\ref{eq:genbound}) are used. Here, $m := 0.5 \overline m$, with $\overline m$ recalculated at each iteration as the maximal value of $m$ that leads to an infeasible problem for (\ref{eq:feasprojrob}). This value is easy to determine via a simple bisection algorithm and, on average, solving 10-20 linear programs with the constraint set (\ref{eq:feasprojrob}). The choice of setting $m$ to half of $\overline m$ is heuristic, and represents a compromise between robustness of the SCFO and convergence speed (as using $\overline m$ directly could lead to very slow progress due to the numerator of (\ref{eq:costfilter2}) being equal to $-\delta_\phi$ consistently).

\end{enumerate}

\begin{figure}
\begin{center}
\includegraphics[width=7cm]{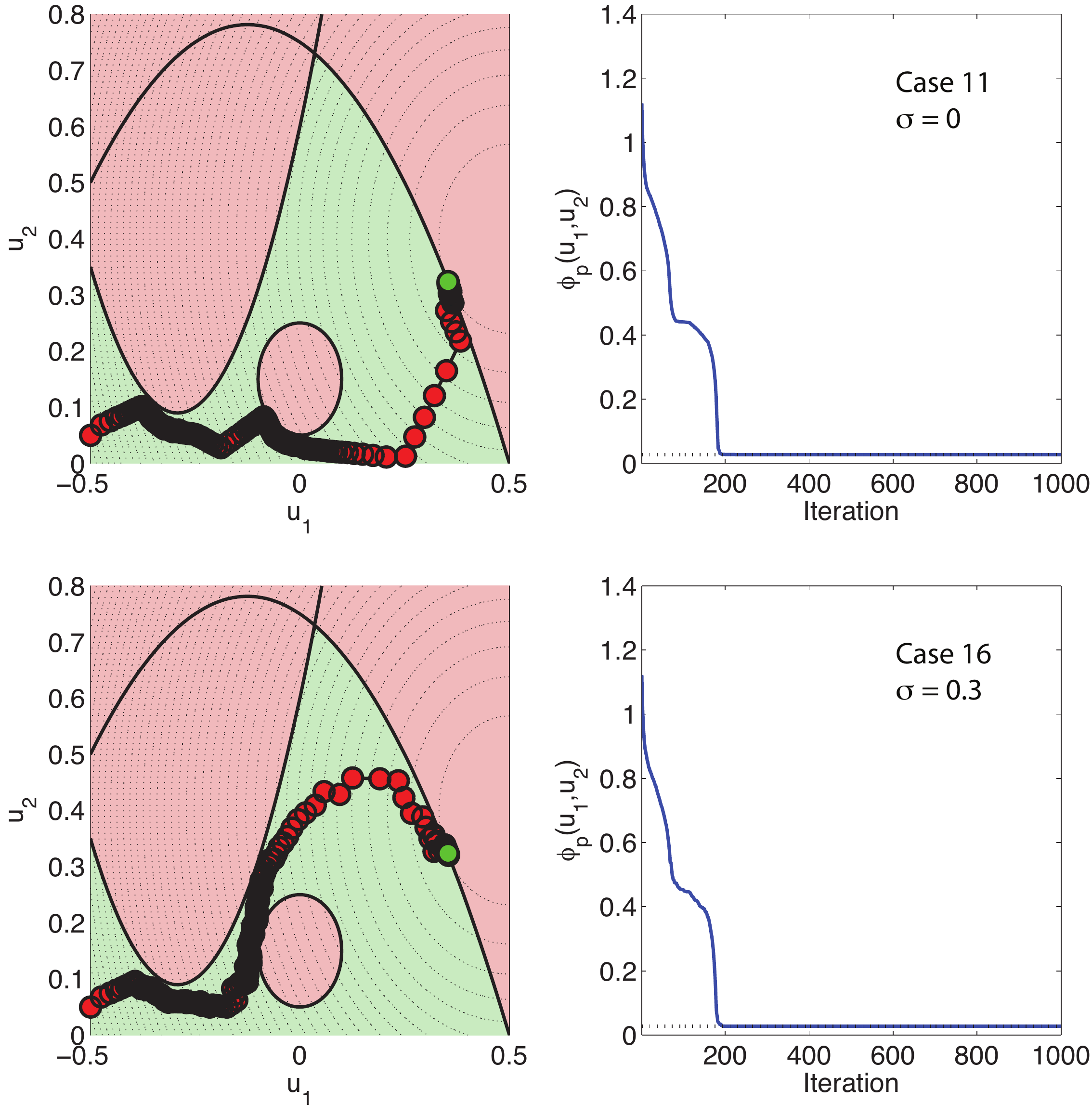}    
\caption{An algorithm may take a different path in the presence of noise and thereby obtain very different performance (here, the performances are nearly identical, despite the presence of noisy gradients in Case 16). As in the previous paper, the red points denote the iterates, the green point denotes the plant optimum, and the thin black lines the contours of the cost function for the left-hand figures. For the figures on the right, the blue line denotes the actual plant cost value at each iteration, while the dotted black line denotes the cost value at the plant optimum.}
\label{fig:P2ex1}
\end{center}
\end{figure}

\begin{figure}
\begin{center}
\includegraphics[width=8cm]{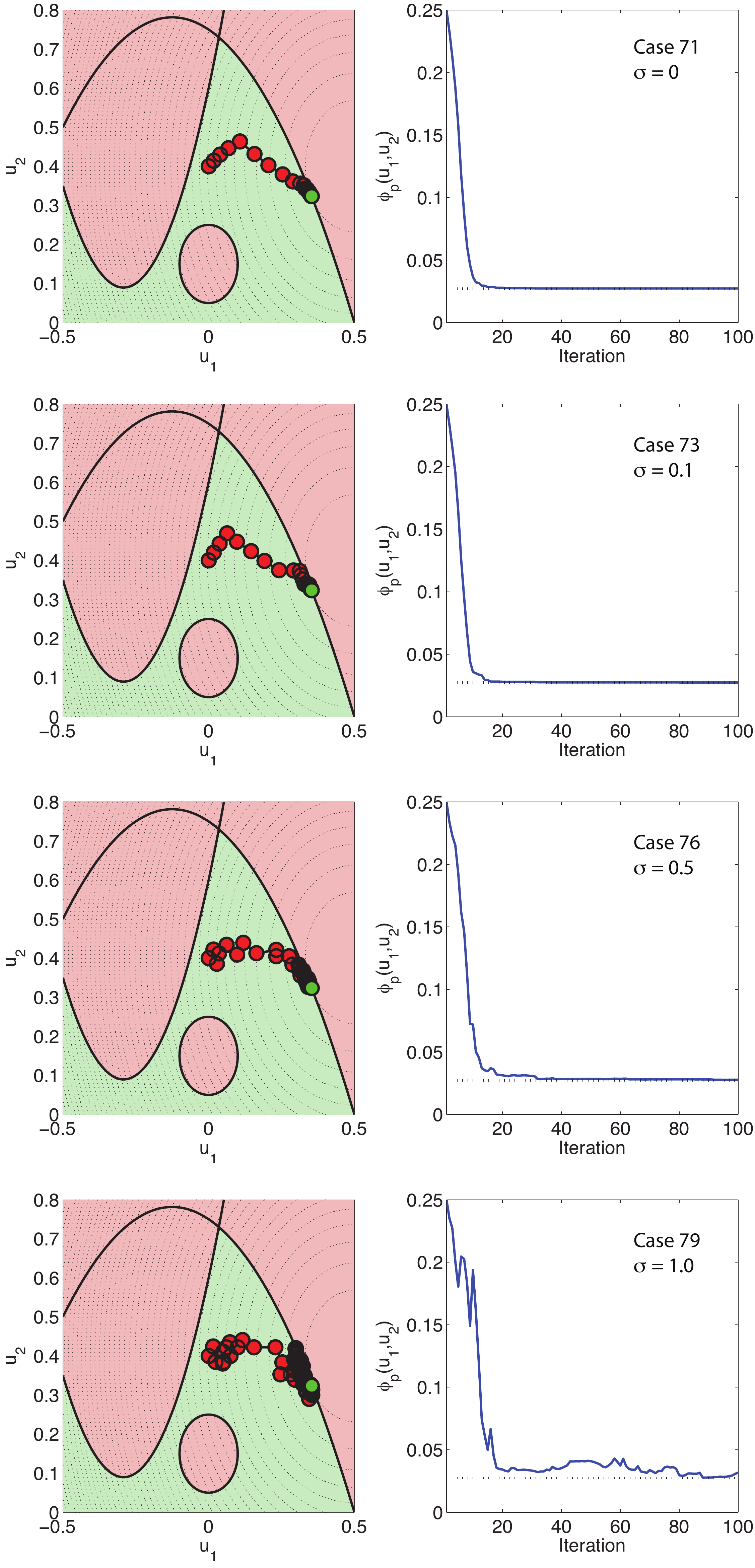}    
\caption{Effect of increasing gradient noise on algorithm performance.}
\label{fig:P2ex2}
\end{center}
\end{figure}

\begin{figure}
\begin{center}
\includegraphics[width=7cm]{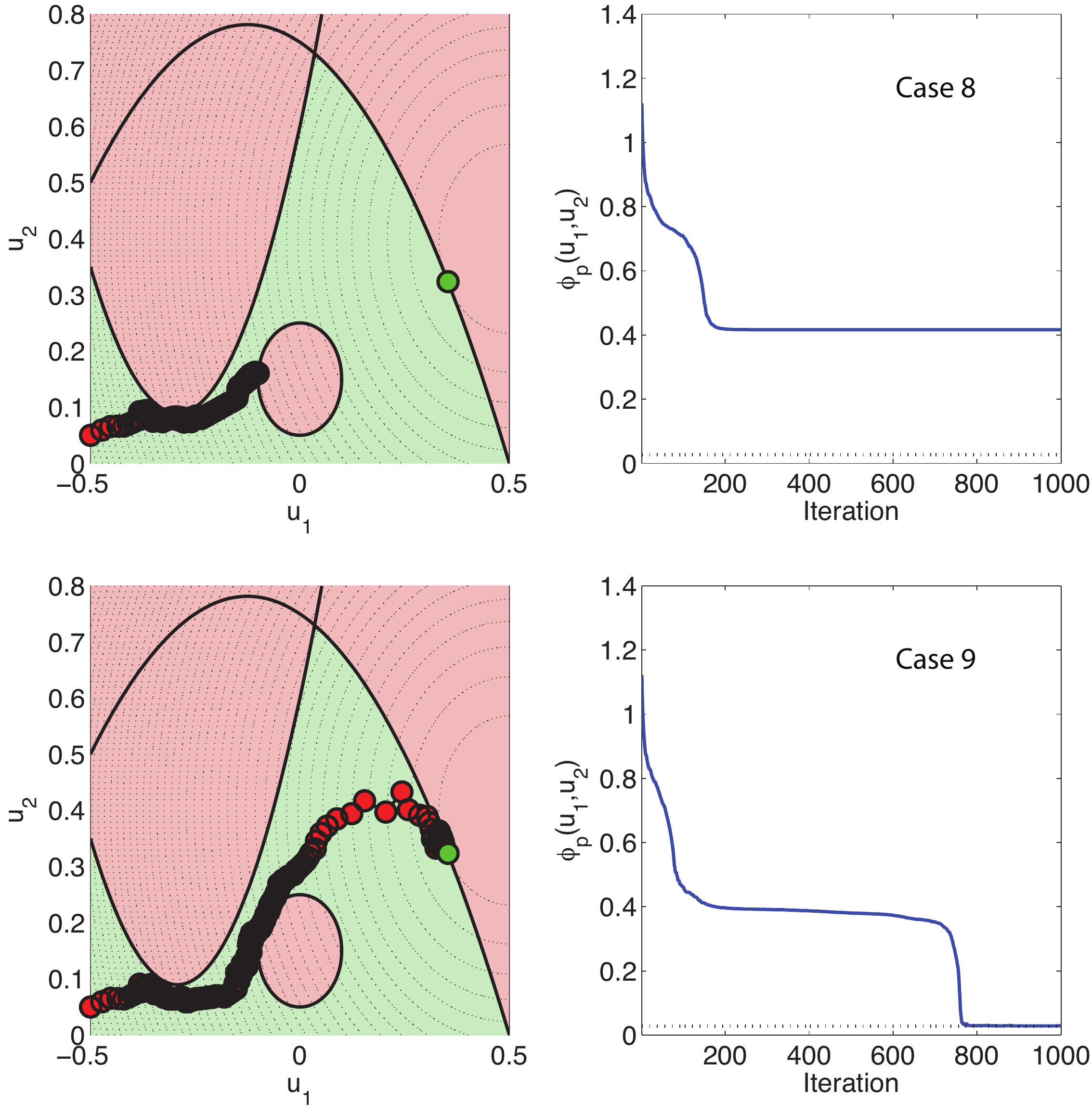}    
\caption{Robust projection (Case 9) may help avoid premature convergence (Case 8).}
\label{fig:P2ex3}
\end{center}
\end{figure}

We present the entire set of results (a total of 100 simulated trials) in Tables \ref{tab:noisygrad1} and \ref{tab:noisygrad2}. While no definitive conclusions are possible, we attempt to highlight the major noticeable trends and to provide insight for them as follows:

\begin{itemize}
\item As expected, using stronger RTO algorithms results in better performance. This may be seen by comparing the performance of the random-step algorithm to the others for both Problems A and B.
\item In general, adding noise to the gradients lowers performance. The sole exception is Case 16, where the algorithm just happens to take an alternate path and obtains relatively better performance as a result (Figure \ref{fig:P2ex1}). A more common effect of the noise is illustrated in Figure \ref{fig:P2ex2}.
\item Implementation III gives consistently better results than Implementation II for Problem A when the noise is moderate ($\sigma = 0.3$) to high ($\sigma = 0.5$). A very likely reason is that the constraints play a major role in this problem, and not avoiding them enough may lead to either premature or slower convergence. As the SCFO and their implementation scheme are designed to encourage a certain back-off from the constraints until approaching them becomes a necessity for cost improvement, it may be expected that robustly enforcing the SCFO will keep the algorithm farther away from the constraints. Not doing so, in turn, would lead to them being approached quicker and thus to slower convergence, which would explain the loss in performance.
\item By contrast, Implementation III only gives consistently better performance for Problem B when the algorithm is very poor (i.e. the random-step algorithm). Otherwise, Implementation II consistently does better. While this may seem puzzling, it may be explained by the simple observation that the SCFO are not crucial for this particular problem and may, in fact, make things worse. Unlike Problem A, the presence of only a single relevant constraint makes robust projections less vital with respect to performance, in that it is easier for the algorithm to navigate the feasible region without running into the constraint, even with corrupted estimates. At the same time, the inherent drawback of the projection used to enforce the SCFO is that it compresses the feasible space of the RTO target based on \emph{local} gradient information, which, when the initial target is very good (as is the case with a good RTO algorithm), may divert this target into a direction that is not as optimal. With a bad algorithm, however, such a diversion does not hurt and may indeed lead to better performance than if no robustness were encouraged, and this is seen with the random-step algorithm.
\item Premature convergence becomes more difficult to avoid when the noise is high and Implementation II is used. However, using Implementation III in this case may help avoid it (see Figure \ref{fig:P2ex3}). Premature convergence is often seen for Problem A when the noise is high.
\item The difficulty in avoiding premature convergence to $g_{p,3}$ (as in Figure \ref{fig:P2ex3}) may be explained by Theorem 1, as this region of the input space contains an (unstable) KKT point at ${\bf u} = [-0.09, 0.11]$, and it is sufficient for the plant gradients corresponding to this point to be included in any $\mathcal{M}_k$ generated at points close by to render the projection infeasible, which in turn leads to reductions in $P$ and thus less robustness. This suggests that unstable KKT points may ``obtain'' a certain stability due to the uncertainty in the gradients.
\item Using Implementation IV leads to somewhat haphazard performance, but this is likely due to the somewhat brute choice of gradient uncertainty, which was, again, taken as $\sigma_\nabla = {\bf 1}$ for all functions. In many cases, it seems to yield optimality losses that are somewhere in between those achieved by Implementations II and III. There are also occasions when it outperforms both, however (see Cases 37, 47, 97, 100). This is, in some sense, encouraging as it shows that a method with significantly fewer assumptions and a fairly brute choice of uncertainty may nevertheless yield acceptable results.

\end{itemize}

\begin{table}
\caption{Performance of different schemes with noisy gradient estimates for Problem A. The following abbreviations are used: IT -- Ideal Target, GD -- Gradient Descent, MA -- Modifier Adaptation, TS -- Two Step, RS -- Random Step, Algo. -- Algorithm, Imp. -- Implementation. Asterisks denote the cases that did not reach the neighborhood of the optimum due to premature convergence to the third constraint.}
\begin{center}
\vspace{-3mm}
\footnotesize{
\begin{tabular}{|p{.6cm}|p{.6cm}|p{.6cm}|p{.6cm}|p{1cm}|}
\hline
Case & Algo. & $\sigma$ & Imp. & $L$   \\ \hline
1 & IT & 0 & I & 73.54 \\ \hline
2 & IT & 0.1 & II & 75.24 \\ \hline 
3 & IT & 0.1 & III & 76.03 \\ \hline 
4 & IT & 0.1 & IV & 75.29 \\ \hline 
5 & IT & 0.3 & II & 118.04 \\ \hline 
6 & IT & 0.3 & III & 82.16 \\ \hline 
7 & IT & 0.3 & IV & 85.06 \\ \hline 
8 & IT & 0.5 & II & 437.69* \\ \hline 
9 & IT & 0.5 & III & 300.44 \\ \hline 
10 & IT & 0.5 & IV & 442.52* \\ \hline \hline
11 & GD & 0 & I & 94.05 \\ \hline
12 & GD & 0.1 & II & 116.61 \\ \hline 
13 & GD & 0.1 & III & 125.42 \\ \hline 
14 & GD & 0.1 & IV & 108.73 \\ \hline 
15 & GD & 0.3 & II & 171.15 \\ \hline 
16 & GD & 0.3 & III & 93.29 \\ \hline 
17 & GD & 0.3 & IV & 124.99 \\ \hline 
18 & GD & 0.5 & II & 431.06* \\ \hline 
19 & GD & 0.5 & III & 434.73* \\ \hline 
20 & GD & 0.5 & IV & 428.58* \\ \hline \hline
21 & MA & 0 & I & 67.70 \\ \hline
22 & MA & 0.1 & II & 67.10 \\ \hline 
23 & MA & 0.1 & III & 68.74 \\ \hline 
24 & MA & 0.1 & IV & 69.08 \\ \hline 
25 & MA & 0.3 & II & 102.00 \\ \hline 
26 & MA & 0.3 & III & 77.34 \\ \hline 
27 & MA & 0.3 & IV & 84.66 \\ \hline 
28 & MA & 0.5 & II & 549.77* \\ \hline 
29 & MA & 0.5 & III & 470.45* \\ \hline 
30 & MA & 0.5 & IV & 433.94* \\ \hline \hline
31 & TS & 0 & I & 76.40 \\ \hline
32 & TS & 0.1 & II & 77.46 \\ \hline 
33 & TS & 0.1 & III & 78.83 \\ \hline 
34 & TS & 0.1 & IV & 75.84 \\ \hline 
35 & TS & 0.3 & II & 102.48 \\ \hline 
36 & TS & 0.3 & III & 94.76 \\ \hline 
37 & TS & 0.3 & IV & 86.42 \\ \hline 
38 & TS & 0.5 & II & 449.62* \\ \hline 
39 & TS & 0.5 & III & 447.36* \\ \hline 
40 & TS & 0.5 & IV & 444.23* \\ \hline \hline 
41 & RS & 0 & I & 82.40 \\ \hline
42 & RS & 0.1 & II & 84.83 \\ \hline 
43 & RS & 0.1 & III & 92.52 \\ \hline 
44 & RS & 0.1 & IV & 79.32 \\ \hline 
45 & RS & 0.3 & II & 124.25 \\ \hline 
46 & RS & 0.3 & III & 102.02 \\ \hline 
47 & RS & 0.3 & IV & 97.83 \\ \hline 
48 & RS & 0.5 & II & 304.69 \\ \hline 
49 & RS & 0.5 & III & 164.62 \\ \hline 
50 & RS & 0.5 & IV & 195.73 \\ \hline 
\end{tabular}
}
\end{center}
\label{tab:noisygrad1}
\end{table}

\begin{table}
\caption{Performance of different schemes with noisy gradient estimates for Problem B.}
\begin{center}
\footnotesize{
\begin{tabular}{|p{.6cm}|p{.6cm}|p{.6cm}|p{.6cm}|p{1cm}|}
\hline
Case & Algo. & $\sigma$ & Imp. & $L$   \\ \hline
51 & IT & 0 & I & 1.12 \\ \hline
52 & IT & 0.1 & II & 1.15 \\ \hline 
53 & IT & 0.1 & III & 1.18 \\ \hline 
54 & IT & 0.1 & IV & 1.18 \\ \hline 
55 & IT & 0.5 & II & 1.26 \\ \hline 
56 & IT & 0.5 & III & 1.62 \\ \hline 
57 & IT & 0.5 & IV & 1.31 \\ \hline 
58 & IT & 1.0 & II & 2.20 \\ \hline 
59 & IT & 1.0 & III & 2.54 \\ \hline 
60 & IT & 1.0 & IV & 2.39 \\ \hline \hline
61 & GD & 0 & I & 1.15 \\ \hline
62 & GD & 0.1 & II & 1.16 \\ \hline 
63 & GD & 0.1 & III & 1.26 \\ \hline 
64 & GD & 0.1 & IV & 1.23 \\ \hline 
65 & GD & 0.5 & II & 1.30 \\ \hline 
66 & GD & 0.5 & III & 1.46 \\ \hline 
67 & GD & 0.5 & IV & 1.19 \\ \hline 
68 & GD & 1.0 & II & 2.00 \\ \hline 
69 & GD & 1.0 & III & 2.45 \\ \hline 
70 & GD & 1.0 & IV & 2.36 \\ \hline \hline
71 & MA & 0 & I & 1.14 \\ \hline
72 & MA & 0.1 & II & 1.16 \\ \hline 
73 & MA & 0.1 & III & 1.19 \\ \hline 
74 & MA & 0.1 & IV & 1.19 \\ \hline 
75 & MA & 0.5 & II & 1.31 \\ \hline 
76 & MA & 0.5 & III & 1.60 \\ \hline 
77 & MA & 0.5 & IV & 1.39 \\ \hline 
78 & MA & 1.0 & II & 2.56 \\ \hline 
79 & MA & 1.0 & III & 2.71 \\ \hline 
80 & MA & 1.0 & IV & 2.48 \\ \hline \hline
81 & TS & 0 & I & 1.19 \\ \hline
82 & TS & 0.1 & II & 1.24 \\ \hline 
83 & TS & 0.1 & III & 1.23 \\ \hline 
84 & TS & 0.1 & IV & 1.20 \\ \hline 
85 & TS & 0.5 & II & 1.39 \\ \hline 
86 & TS & 0.5 & III & 1.75 \\ \hline 
87 & TS & 0.5 & IV & 1.36 \\ \hline 
88 & TS & 1.0 & II & 2.83 \\ \hline 
89 & TS & 1.0 & III & 2.58 \\ \hline 
90 & TS & 1.0 & IV & 2.64 \\ \hline \hline 
91 & RS & 0 & I & 1.51 \\ \hline
92 & RS & 0.1 & II & 1.52 \\ \hline 
93 & RS & 0.1 & III & 1.44 \\ \hline 
94 & RS & 0.1 & IV & 1.23 \\ \hline 
95 & RS & 0.5 & II & 1.76 \\ \hline 
96 & RS & 0.5 & III & 1.49 \\ \hline 
97 & RS & 0.5 & IV & 1.27 \\ \hline 
98 & RS & 1.0 & II & 2.81 \\ \hline 
99 & RS & 1.0 & III & 2.03 \\ \hline 
100 & RS & 1.0 & IV & 1.65 \\ \hline 
\end{tabular}
}
\end{center}
\label{tab:noisygrad2}
\end{table}

\section{Inaccurate Measurement or Estimation of the Constraints}

We recall that Conditions (\ref{eq:suffcond2}) and (\ref{eq:suffcond3}) depend on knowing the values of $g_{p,j}({\bf u}_k)$. In practice, however, the constraint value $g_{p,j} ({\bf{u}}_k)$ may not be measured directly. It could, for example, be estimated from (measured) outputs. Additionally, even when measured directly via a sensor, there is likely to be noise and measurement errors that prevent us from having accurate knowledge of the real constraint value. When the values are inaccurate, it follows that the upper bound values on the filter gains may be higher than they truly are, which in turn jeopardizes the feasibility guarantees. There is also the chance of premature convergence as the Boolean trigger in (\ref{eq:suffcond2}) may not be activated when it should be.

Here, we proceed to derive the robust version of these conditions by proposing different ways to bound the value of $g_{p,j}({\bf u}_k)$ at each iteration, and demonstrate their effectiveness with a number of simulated trials.

\subsection{Bounding Constraint Values}

Let a constraint value for the iteration $k$ be bounded from above with some satisfactory confidence as:

\begin{equation}\label{eq:boundedcon}
g_{p,j}({\bf{u}}_k) \leq \overline g_{p,j}({\bf{u}}_k).
\end{equation}

Since:

\begin{equation}\label{eq:gainlow}
\frac{{- \overline g_{p,j}({\bf{u}}_{k } )}}{\displaystyle \sum\limits_{i = 1}^{n_u} {\kappa_{ji} | u^*_{k+1,i} - u_{k,i} | }} \leq \frac{{- g_{p,j}({\bf{u}}_{k } )}}{\displaystyle \sum\limits_{i = 1}^{n_u} {\kappa_{ji} | u^*_{k+1,i} - u_{k,i} | }} ,
\end{equation}

\noindent it follows that:

\begin{equation}\label{eq:gainupperrob}
\begin{array}{l}
K_{k} \le \mathop{\min}\limits_{j = 1,...,n_g} \left[\displaystyle \frac{{- \overline g_{p,j}({\bf{u}}_{k } )}}{\displaystyle \sum\limits_{i = 1}^{n_u} {\kappa_{ji} | u^*_{k+1,i} - u_{k,i} | }} \right] \vspace{2mm}\\
\hspace{10mm}\Rightarrow K_{k} \le \mathop{\min}\limits_{j = 1,...,n_g} \left[ \displaystyle \frac{{- g_{p,j}({\bf{u}}_{k } )}}{\displaystyle \sum\limits_{i = 1}^{n_u} {\kappa_{ji} | u^*_{k+1,i} - u_{k,i} | }} \right]
\end{array},
\end{equation}

\noindent with the left-hand side being the \emph{robust feasibility condition}.

We note that (\ref{eq:gainupperrob}), while conceptually robust, may suffer considerable practical drawbacks if the upper bound on the true constraint value is very conservative. As such, we now propose several ways in which this bound may be chosen and how, if necessary, it may be refined so as to yield the desired performance with increasing iterations.

\begin{thm}{\bf (Robust Upper Bounds on Constraint Values)}

Any of the following choices for $\overline g_{p,j}({\bf{u}}_k)$ satisfy (\ref{eq:boundedcon}):

\begin{enumerate}[(i)]
\item $\overline g_{p,j}({\bf{u}}_k) = 0$.
\item $\overline g_{p,j} ({\bf{u}}_k) = \hat g_{p,j} ({\bf{u}}_k) - \underline w_j$, where $\hat g_{p,j} ({\bf{u}}_k)$ is the measured or estimated value of the constraint at iteration $k$, and $w_{k,j}$ is an additive noise/error term, with the lower bound $\underline w_j \leq w_{k,j}, \; \forall k$.
\item $\overline g_{p,j} ({\bf{u}}_k) = \frac{1}{n} \Sigma_{i = k-n+1}^{k}\; \hat g_{p,j}({\bf{u}}_{i}) - \underline w_j (n)$, where ${\bf u}_k = {\bf u}_{k-1} = ... = {\bf u}_{k-n+1}$ and $\underline w_j (n) \leq  \frac{1}{n} \Sigma_{i = k-n+1}^{k} w_{i,j}$ is a lower bound on the mean of $n$ noise/error terms.
\item $\overline g_{p,j} ({\bf{u}}_k) = \overline g_{p,j} ({\bf{u}}) + \displaystyle \sum\limits_{i = 1}^{n_u} {\kappa_{ji} | u_{k,i} - u_{i} | }, \; \forall {\bf u} \in \mathcal{I}$.
\end{enumerate}

\end{thm}

\begin{pf}

\begin{enumerate}[(i)]
\item Follows from the robust guarantee of feasibility due to (\ref{eq:gainupperrob}).
\item Follows from the definition of $\hat g_{p,j} ({\bf{u}}_k)$ and the lower bound on $w_{k,j}$:

\begin{equation}\label{eq:th2A}
\begin{array}{l}
\hat g_{p,j} ({\bf{u}}_k) = g_{p,j} ({\bf{u}}_k) + w_{k,j}  \\
\Rightarrow g_{p,j} ({\bf{u}}_k) = \hat g_{p,j} ({\bf{u}}_k) - w_{k,j} \leq \hat g_{p,j} ({\bf{u}}_k) - \underline w_j
\end{array}.
\end{equation}

\item Follows a similar derivation to (ii):

\begin{equation}\label{eq:th2B}
\begin{array}{l}
\hat g_{p,j} ({\bf{u}}_{k-i+1}) = g_{p,j} ({\bf{u}}_k) + w_{k-i+1,j}, \;\;\; i = 1,...,n  \\
\Rightarrow n g_{p,j} ({\bf{u}}_k) = \displaystyle \sum_{i = k-n+1}^{k}\; \hat g_{p,j}({\bf{u}}_{i}) - \displaystyle \sum_{i = k-n+1}^{k} w_{i,j}  \\
\Rightarrow g_{p,j} ({\bf{u}}_k) = \displaystyle \frac{1}{n} \sum_{i = k-n+1}^{k}\; \hat g_{p,j}({\bf{u}}_{i}) - \displaystyle \frac{1}{n} \sum_{i = k-n+1}^{k} w_{i,j}  \\
\Rightarrow g_{p,j} ({\bf{u}}_k) \leq \displaystyle \frac{1}{n} \sum_{i = k-n+1}^{k}\; \hat g_{p,j}({\bf{u}}_{i}) - \underline w_j(n)
\end{array}.
\end{equation}

\item Follows directly from the Lipschitz bound (\ref{eq:linbound}). \qed
\end{enumerate} 

\end{pf}

As all of these four bounds are valid, we may simply take their minimum as $\overline g_{p,j}({\bf{u}}_k)$ in implementation. We make the following remarks:

\begin{itemize}
\item Bound (i) is not useful as it essentially results in $K_k = 0$. However, it is important since it guarantees that feasibility is kept in the worst case, when Bounds (ii)-(iv) cannot guarantee $\overline g_{p,j}({\bf{u}}_k) < 0$.
\item The nature of Bounds (ii) and (iii) is largely stochastic if one assumes that the noise comes from some probability distribution. In many cases, we may expect that $\overline g_{p,j}({\bf{u}}_k) \rightarrow g_{p,j}({\bf{u}}_k)$ as $n \rightarrow \infty$ in Bound (iii). Bound (iii) would naturally come into play whenever Bound (i) is used for consecutive iterations. For obtaining the lower bounds on the noise, we note that $\underline w_j$ may be obtained directly from the probability distribution of the noise, while $\underline w_j (n)$ may be approximated using Monte Carlo sampling.
\item Bound (iv) is innately deterministic, but becomes stochastic when Bounds (ii) and (iii) from previous measurements are used.
\end{itemize}
 
Also important is the Boolean nature of the projection condition:

\begin{equation}\label{eq:boolean}
\forall j: g_{p,j}({\bf{u}}_k) \geq -\epsilon_j,
\end{equation}

\noindent which will be corrupted if $\hat g_{p,j}({\bf{u}}_k)$ is used instead. Although this condition is met ``on average'' in some cases (i.e. when the noise or error is zero-mean), we propose to make it robust by implementing the following instead:

\begin{equation}\label{eq:booleanrob}
\forall j: \overline g_{p,j}({\bf{u}}_k) \geq -\epsilon_j,
\end{equation}

\noindent as this guarantees that a local descent in the constraint is always enforced when the true value is $\epsilon$-active, and thereby precludes the possibility of premature convergence due to inaccurate constraint values. 

\subsection{Examples}

We consider the same two problems as before, but this time corrupt the constraint measurements as in (\ref{eq:th2A}), with $w_{k,j} \sim \mathcal{N}(0,(\sigma_g \overline \epsilon_{j})^2)$, which essentially means that the noise is somehow proportional (via the constant $\sigma_g$, which is varied for test purposes) to the maximum absolute value, or the range, of the constraint on the relevant input space. The lower bound $\underline w_j$ is chosen as $-3\sigma_g \overline \epsilon_{j}$, which corresponds to 99.85\% confidence, and the standard variance reduction law for a normally distributed noise, $\underline w_j (n) = -3\sigma_g \overline \epsilon_{j} / \sqrt{n}$, is used.

\begin{figure}
\begin{center}
\includegraphics[width=8cm]{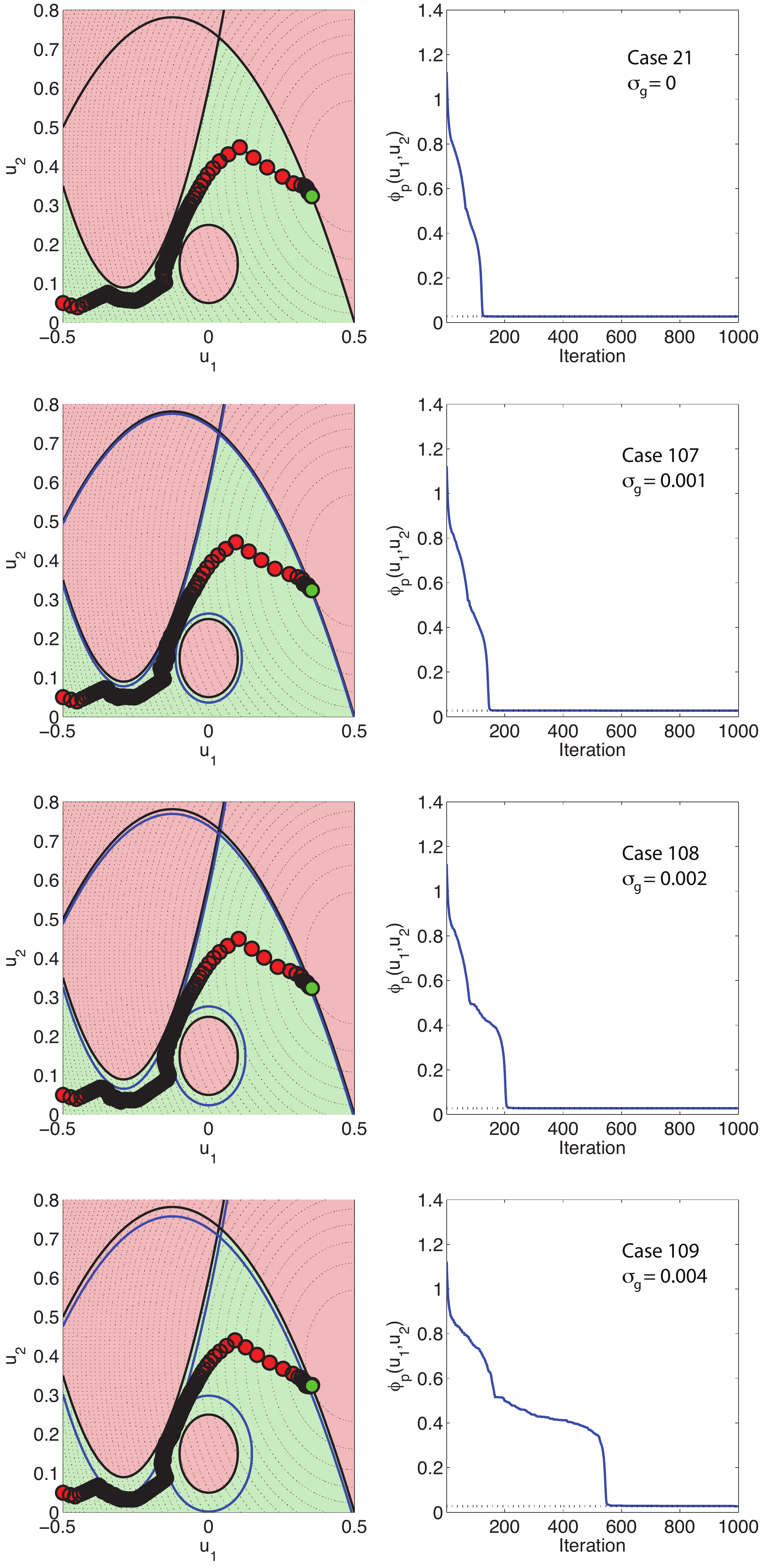}    
\caption{RTO performance for Problem A with increasing noise levels in the constraint measurements. The magnitude of the noise is illustrated by giving, in blue, the contours of the constraints with back-offs equal to $-\underline w_j$.}
\label{fig:P2ex6}
\end{center}
\end{figure}

\begin{figure}
\begin{center}
\includegraphics[width=8cm]{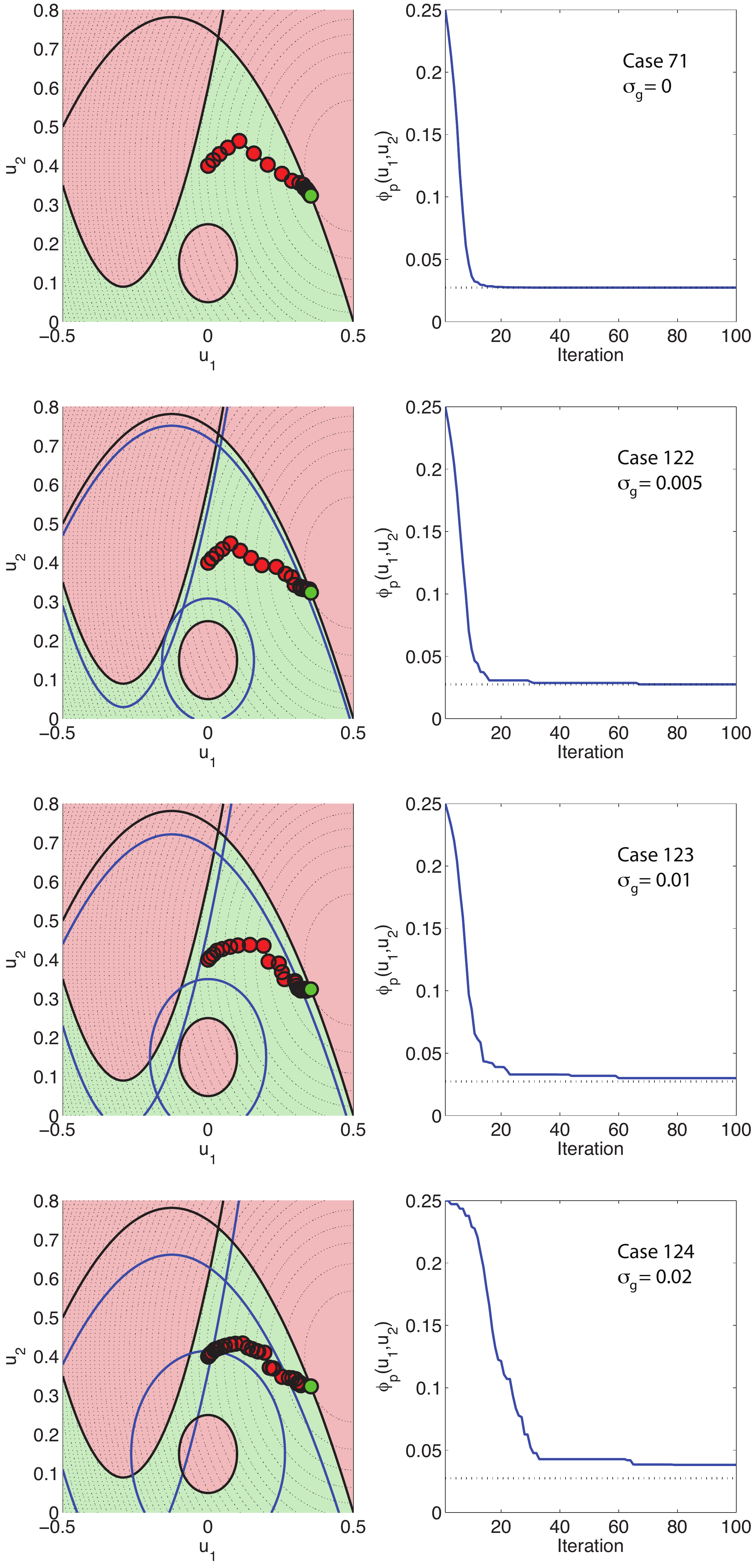}    
\caption{RTO performance for Problem B with increasing noise levels in the constraint measurements.}
\label{fig:P2ex7}
\end{center}
\end{figure}

\begin{figure}
\begin{center}
\includegraphics[width=8cm]{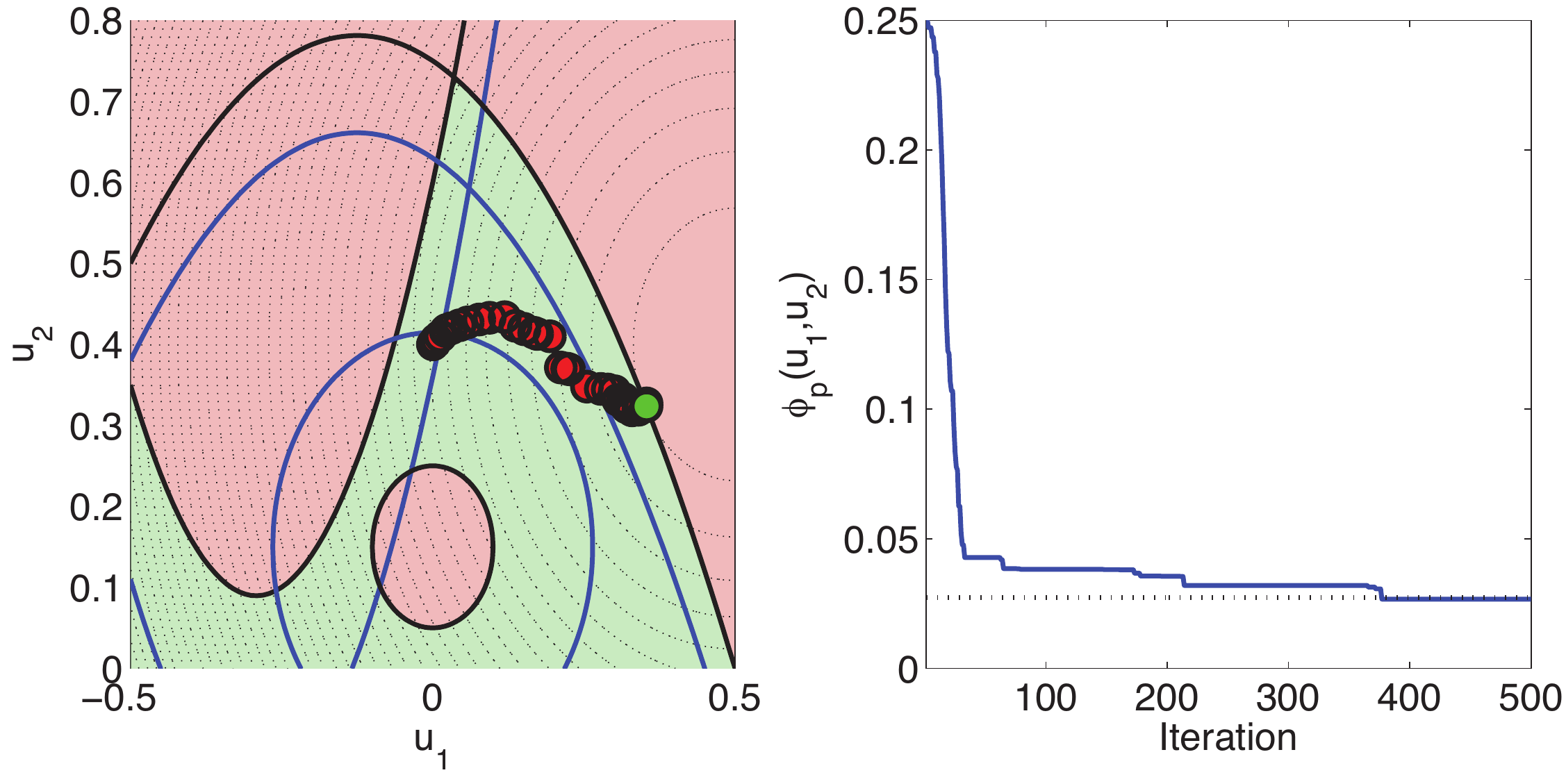}    
\caption{A longer simulation of Case 124 showing that the offset from the optimum may be eliminated as the noise is filtered out with increasing iterations.}
\label{fig:P2ex9}
\end{center}
\end{figure}

The results of the numerical trials are presented in Tables \ref{tab:noisycon1} and \ref{tab:noisycon2}, and are in line with what would be expected, as augmenting the noise in all cases leads to consistently poorer performance. We illustrate this through the examples in Figures \ref{fig:P2ex6} and \ref{fig:P2ex7}. To show that the algorithm is able to approach the plant optimum asymptotically, even when the noise is quite large, we re-run the problem of Case 124 for 500 iterations instead of 100, which shows that the offset from the optimum is gradually reduced as the noise is filtered out (Figure \ref{fig:P2ex9}). Although hard to see in this figure, we note that the final set of iterates actually \emph{violates} the constraint, if only slightly. This is, however, due to the probabilistic nature of Bounds (ii) and (iii) and could be remedied by choosing an even higher confidence level.

\begin{table}
\caption{Performance of different schemes with noisy constraint estimates for Problem A.}
\begin{center}
\footnotesize{
\begin{tabular}{|p{.6cm}|p{.6cm}|p{.6cm}|p{1cm}|}
\hline
Case & Algo. & $\sigma_g$ & $L$   \\ \hline
1 & IT & 0 & 73.54 \\ \hline
101 & IT & 0.001 & 84.54 \\ \hline 
102 & IT & 0.002 & 112.48 \\ \hline 
103 & IT & 0.004 & 304.14 \\ \hline \hline
11 & GD & 0 & 94.05 \\ \hline
104 & GD & 0.001 & 108.17 \\ \hline 
105 & GD & 0.002 & 154.69 \\ \hline 
106 & GD & 0.004 & 476.42 \\ \hline \hline 
21 & MA & 0 & 67.70 \\ \hline
107 & MA & 0.001 & 80.06 \\ \hline 
108 & MA & 0.002 & 105.81 \\ \hline 
109 & MA & 0.004 & 271.59 \\ \hline  \hline
31 & TS & 0 & 76.40 \\ \hline
110 & TS & 0.001 & 105.45 \\ \hline 
111 & TS & 0.002 & 121.75 \\ \hline 
112 & TS & 0.004 & 315.49 \\ \hline  \hline 
41 & RS & 0 & 82.40 \\ \hline
113 & RS & 0.001 & 117.59 \\ \hline 
114 & RS & 0.002 & 127.77 \\ \hline 
115 & RS & 0.004 & 393.02 \\ \hline  
\end{tabular}
}
\end{center}
\label{tab:noisycon1}
\end{table}

\begin{table}
\caption{Performance of different schemes with noisy constraint estimates for Problem B.}
\begin{center}
\footnotesize{
\begin{tabular}{|p{.6cm}|p{.6cm}|p{.6cm}|p{1cm}|}
\hline
Case & Algo. & $\sigma_g$ & $L$   \\ \hline
51 & IT & 0 & 1.12 \\ \hline
116 & IT & 0.005 & 1.60 \\ \hline 
117 & IT & 0.010 & 1.91 \\ \hline 
118 & IT & 0.020 & 5.02 \\ \hline \hline
61 & GD & 0 & 1.15 \\ \hline
119 & GD & 0.005 & 1.63 \\ \hline 
120 & GD & 0.010 & 2.23 \\ \hline 
121 & GD & 0.020 & 6.47 \\ \hline \hline 
71 & MA & 0 & 1.14 \\ \hline
122 & MA & 0.005 & 1.48 \\ \hline 
123 & MA & 0.010 & 2.04 \\ \hline 
124 & MA & 0.020 & 5.09 \\ \hline  \hline
81 & TS & 0 & 1.19 \\ \hline
125 & TS & 0.005 & 1.62 \\ \hline 
126 & TS & 0.010 & 2.20 \\ \hline 
127 & TS & 0.020 & 6.45 \\ \hline  \hline 
91 & RS & 0 & 1.51 \\ \hline
128 & RS & 0.005 & 1.95 \\ \hline 
129 & RS & 0.010 & 2.93 \\ \hline 
130 & RS & 0.020 & 7.37 \\ \hline  
\end{tabular}
}
\end{center}
\label{tab:noisycon2}
\end{table}

\section{Lipschitz Constants and the Quadratic Upper Bound}

It is clear that both the Lipschitz constants for the constraints and the quadratic upper bound on the cost, as defined in (\ref{eq:linbound}) and (\ref{eq:Qbound}), play vital roles in enforcing the SCFO by acting to determine an input filter gain that both preserves feasibility and yields a step that decreases the cost at the following iteration. Additionally, the Lipschitz constants may play other roles as well, such as serving as the ultimate lower and upper bounds on the gradient estimate, as in (\ref{eq:lip0}), and helping define an upper bound on a constraint function value as proposed in Theorem 2. However, depending on the application, the choice of a given Lipschitz constant or, further yet, the values of the quadratic upper bound, may not be trivial.

In this section, we treat these issues in some detail, starting with the idea of a Lipschitz constant and demonstrating that, despite its somewhat mathematical origin, it has a physical significance that may become apparent quite readily in some applications. We then extend to the cases where this is not so, and propose other means to define the constants, or to relax the resulting upper bounds by using all of the data obtained by the RTO algorithm during operation. We finish by carrying out a similar discussion for the quadratic upper bound, and look at the cases where it may or may not be crucial. Examples are given throughout the section to supplement the presented ideas where necessary.

\subsection{The Physical Meaning of Lipschitz Constants}

The vast majority of mathematical works present the (strict) Lipschitz constant as a single constant, $\kappa_j \in \mathbb{R}$, that meets the following criterion:

\begin{equation}\label{eq:lipstan}
\begin{array}{l}
| g_{p,j} ({\bf{u}}_{k+1}) - g_{p,j} ({\bf{u}}_k) | < \kappa_j \| {\bf u}_{k+1} - {\bf u}_{k}  \| \vspace{1mm}\\
\hspace{10mm}\forall {\bf u}_{k},{\bf u}_{k+1}  \in \mathcal{I} \setminus \left\{ {\bf u}_{k}, {\bf u}_{k+1} : {\bf u}_{k} = {\bf u}_{k+1} \right\} 
\end{array},
\end{equation}

\noindent for some choice of norm. While useful for many conceptual reasons, such as proving that certain optimization algorithms converge in finite time \citep[e.g.][Section 2.5]{Fletcher:87}, or even concrete algorithms, such as the Lipschitz optimization branch of global optimization \citep[Ch. 5]{Horst:95}, the Lipschitz constant as defined in (\ref{eq:lipstan}) remains a very abstract thing -- when $g_{p,j} ({\bf{u}})$ is a real plant constraint, it is difficult to find a physical and understandable interpretation to $\kappa_j$ as defined in (\ref{eq:lipstan}), apart from saying that it is the ``worst-case relative change in the constraint for a perturbation in the inputs with respect to some norm''. As a result, even an experienced plant operator or engineer -- one who is well familiar with the system at hand -- may find it difficult to provide a value of $\kappa_j$ that is not overly conservative. Additionally, while some work has been done on estimating Lipschitz constants in the numerical context \citep{Meewella:88,Hansen:92,Wood:96}, these are often limited to low-dimensional cases and would be difficult to apply in RTO.

This is a major reason for our choosing to work with univariate Lipschitz constants as defined by (\ref{eq:linbound}) and (\ref{eq:lip0}). As stated in the latter, the practical meaning of a given Lipschitz constant when defined this way is quite clear -- it is simply the maximal sensitivity of a single constraint to a single input, or, in layman's terms, it is ``the largest relative increase that one can observe in the constraint value when perturbing a single given input''. Such a definition, we believe, has a much more concrete interpretation, as may be illustrated by the following example.

\vspace{2mm}
\noindent{\it{Example -- Lipschitz Constants for a CSTR with Two Reactions}} 
\vspace{2mm}

\begin{figure}
\begin{center}
\includegraphics[width=8cm]{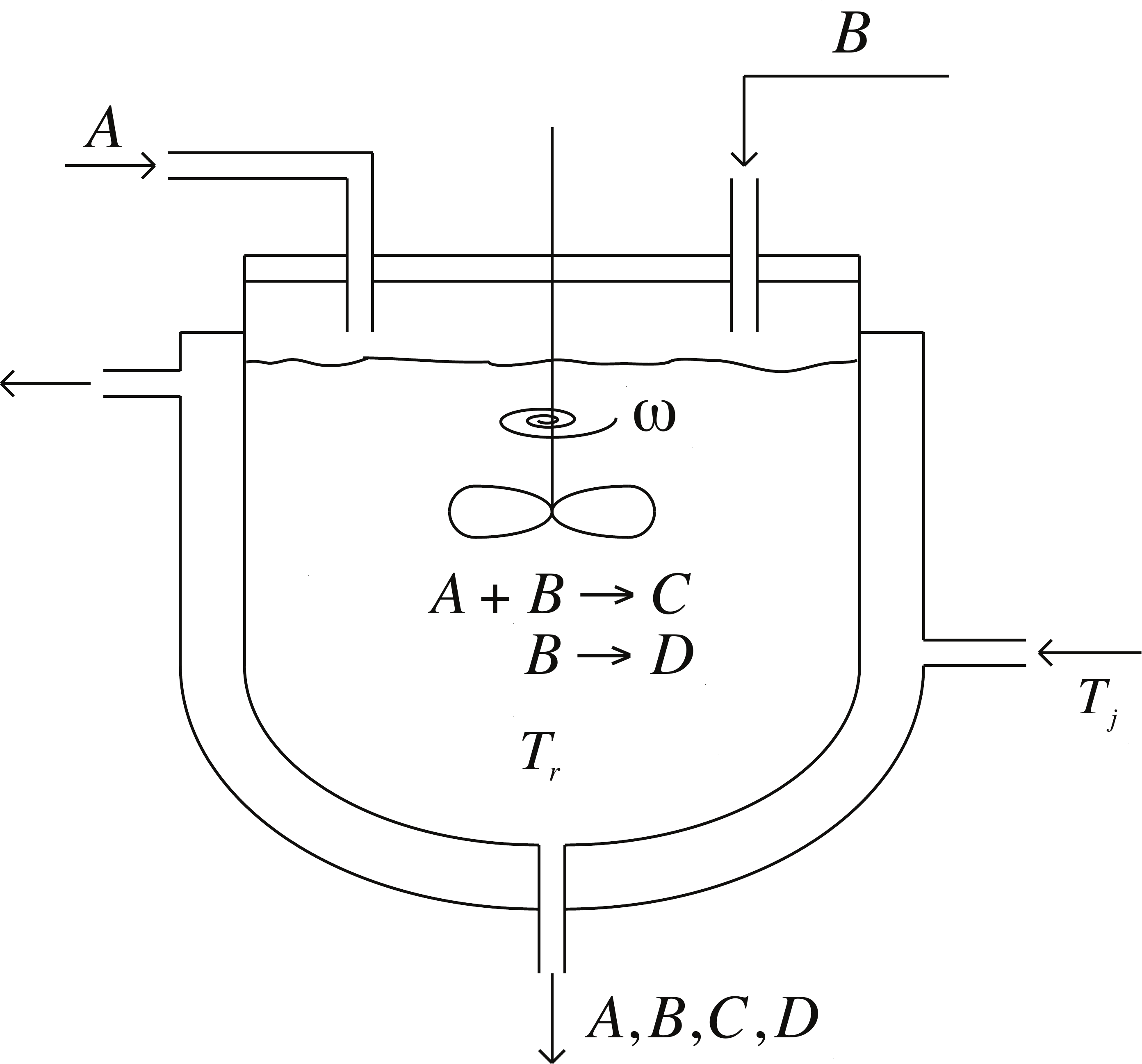}    
\caption{The schematic of the CSTR setup for Problem (\ref{eq:CSTRopt}).}
\label{fig:ex1}
\end{center}
\end{figure}

Consider a continuous-stirred tank reactor (CSTR) that is fed with reactants $A$ and $B$ and has an outlet stream of both reactants and products $A$, $B$, $C$, and $D$, with the reactor temperature, $T_{r}$, controlled by manipulating the temperature of the inlet stream to the surrounding jacket, $T_j$. An impeller mixes the contents with the velocity $\omega$ (Figure \ref{fig:ex1}). Two reactions, $A + B \rightarrow C$ and $B \rightarrow D$, take place in the reactor and have kinetics that vary with both the temperature and the mixing rate. 

Suppose now that the feed rates of $A$ and $B$ are fixed, and the goal is to optimize the steady-state yield, $\gamma$, of $C$. The manipulated variables are chosen as the temperature of the jacket inlet and the velocity of the impeller. Additionally, there is a maximum-temperature constraint on $T_{r}$, past which it is believed that the reactions produce very undesirable results. Product $D$ is also considered as undesirable and the number of moles of $D$ that are produced, $n_D$, is required to be kept below a certain threshold. Adding the fact that both the jacket inlet temperature and the velocity of the impeller are to be kept within certain operating bounds, this translates into the following RTO problem:

\begin{equation}\label{eq:CSTRopt}
\begin{array}{l}
\mathop {{\rm{maximize}}}\limits_{T_j,\omega}\hspace{4mm}\gamma (T_j,\omega)  \\
{\rm{subject}}\hspace{1mm}{\rm{to}}\hspace{3mm}T_{r}(T_j,\omega) \leq T^U_{r} \vspace{1mm}\\
\hspace{18mm}n_D (T_j,\omega) \leq n ^U_{D} \vspace{1mm} \\
\hspace{18mm}T^L_j \leq T_j \leq T^U_j,\hspace{2mm} \omega^L \leq \omega \leq \omega^U
\end{array}\hspace{3mm}.
\end{equation}

This problem has two input variables, $T_j$ and $\omega$, and two uncertain constraint functions, $T_{r}(T_j,\omega)$ and $n_D (T_j,\omega)$, thereby resulting in four Lipschitz constants: $\kappa_{T_{r},T_j}$, $\kappa_{T_{r},\omega}$, $\kappa_{n_D,T_j}$, and $\kappa_{n_D,\omega}$.

Suppose that this is a process where the enthalpies of reaction are close to 0, i.e. the reactions contribute very little to the temperature change in the CSTR. In this case, $\kappa_{T_{r},T_j}$ may be set to 1 -- e.g. increasing the jacket inlet temperature by 10\degree C cannot increase the reactor temperature by more than 10\degree C.

For $\kappa_{T_{r},\omega}$, we may rely on our experience and past observations that changing impeller velocity has negligible effect on the temperature at steady state. We could therefore set $\kappa_{T_{r},\omega}$ to 0, or to a very small value. For these two Lipschitz constants, the choice is relatively simple.

$\kappa_{n_D,T_j}$ and $\kappa_{n_D,\omega}$ require some work, however, since the exact relations between the number of moles of $D$ produced and the jacket inlet temperature and impeller velocity are more complicated than the heat transfer case above. Here, we may fall back on modeling and simulations to compute an estimate of the worst-case sensitivities, and then, for example, take a multiple of these estimates as a further safety precaution. If no model is available, we may fall back on past experimental data, if not for precisely this process then for a similar one, in hopes of obtaining a reasonable estimate. If this is not possible, then an overly conservative large value may be chosen, this latter corresponding to the least desirable scenario.

\vspace{3mm}

The above example serves to show a hierarchy by which one may proceed to choose the Lipschitz constants for a given process:

\begin{enumerate}[i.]
\item Use the physical laws of the process to choose the Lipschitz constants based on natural bounds that cannot be exceeded.
\item Numerically calculate the maximal sensitivities of a model of the process (e.g. a first-principles model). Then, augment these values by a safety factor that depends on the lack of confidence in the model.
\item Study past experimental input-output data to get a sense for what the sensitivities usually are and how large they may become. Then, augment these values by a safety factor that depends on how valid and global these experimental data are.
\item Use overly large, but safe, values for the Lipschitz constants.
\end{enumerate}

In the discussion that follows, we go through some special cases of these methods, all of which are used to choose the Lipschitz constants \emph{prior} to optimization, and then extend on them by showing how data obtained during the RTO iterations may be used to further improve on the bound of (\ref{eq:linbound}).

\subsection{Obtaining and Refining Lipschitz Bounds}

\subsubsection{Exploiting Direction}

In all of the discussion in both this and the companion paper, we have chosen to use the same Lipschitz constant for both the lower and upper bounds on a given derivative, as in (\ref{eq:lip0}). However, there is no reason to limit ourselves to such a definition, and it is generally advised that the distinction between the lower and upper Lipschitz constants be made:

\begin{equation}\label{eq:liptwo}
\underline \kappa_{ji} < \frac{\partial g_{p,j}}{\partial u_i} \Big |_{\bf{u}} < \overline \kappa_{ji}, \;\; \forall {\bf u} \in \mathcal{I},
\end{equation}

\noindent as doing so has the advantage of less conservatism.

As an example, we may return to the CSTR problem and consider the effect of $T_j$ on $T_{r}$, noting now that it is very unlikely that decreasing $T_j$ will ever increase $T_{r}$, which implies that the sensitivity is always positive. We could thereby set $\underline \kappa_{T_{r},T_j} = 0$ and $\overline \kappa_{T_{r},T_j} = 1$. 

Following the same logic as in (\ref{eq:der2}), we may define the refinement analytically:

\begin{equation}\label{eq:lipthree}
\displaystyle \sum\limits_{i = 1}^{n_u} {\kappa_{ji} |  u^*_{k+1,i} - u_{k,i} | } \geq \displaystyle \sum\limits_{i = 1}^{n_u} {\mathop {\max} \left( \underline \kappa_{ji} ( u^*_{k+1,i} - u_{k,i} ), \overline \kappa_{ji} ( u^*_{k+1,i} - u_{k,i} ) \right) },
\end{equation}

\noindent which, in the CSTR example, proves useful as it does not assume a growth in the constraint value regardless of direction (as is done with the absolute value on the left-hand side). For simplicity, we will continue to use the left-hand-side form for the remainder of this document, but emphasize that the right-hand side can, and should, be substituted into any expressions employing the Lipschitz bound.

\subsubsection{Exploiting Sparsity in the Jacobian}

The Jacobian matrix is an $n_g \times n_u$ matrix of constraint-to-input sensitivities, each of which has a corresponding Lipschitz constant. In practice, it may often occur that this matrix is not full and that some of its elements may be 0 due to a complete absence of a relationship between a given constraint and a given input. When such cases can be confirmed with certainty in applications (from natural physical laws or the setup of the problem), the corresponding Lipschitz constants would be 0 as well.

\subsubsection{Exploiting Partial Parametric Uncertainty}

In many applications, a model of the constraints may be available. While it is generally risky to assume that all of the errors in such models may be described fully by the uncertain model parameters, it may nevertheless occur that certain parts of the constraints are fairly well-modeled and may be assumed to err in the parameters only. If this is the case, then the choice of a given Lipschitz constant may be made by solving a numerical optimization problem.

For example, consider a two-input RTO problem where one of the constraints is described by the model:

\begin{equation}\label{eq:paruncertainty}
g_j ({\bf u},\theta) = 2u_1 + \theta u_2^2,
\end{equation}

\noindent and where the relevant input space is defined by $u_1,u_2 \in [0,1]$ and the unknown parameter $\theta$ is assumed to vary between 1 and 5. If we know that the term $\theta u_2^2$ is the only means by which $u_2$ enters into this constraint and are sure that the parameter $\theta$ captures the modeling error fully, then we may calculate the Lipschitz constants $\underline \kappa_{j2}$ and $\overline \kappa_{j2}$ by simply evaluating the following:

\begin{equation}\label{eq:parlip}
\displaystyle \mathop{\inf}\limits_{{\scriptsize \begin{array}{c}\theta \in [1,5]\\ u_1,u_2 \in [0,1] \end{array}}} \frac{\partial g_j ({\bf u},\theta)}{\partial u_2}, \;\;\; \displaystyle \mathop{\sup}\limits_{{\scriptsize \begin{array}{c}\theta \in [1,5]\\ u_1,u_2 \in [0,1] \end{array}}} \frac{\partial g_j ({\bf u},\theta)}{\partial u_2},
\end{equation}

\noindent which, in this case, gives $\underline \kappa_{j2} = 0$ and $\overline \kappa_{j2} = 10$. These values can be further refined online if we only consider a local subspace of the relevant input space, defined by the current iterate and the projected RTO target:

\begin{equation}\label{eq:localspace}
\mathcal{I}_k = \{ {\bf u} : \mathop {\min} \left( u_{k,i}, u_{k+1,i}^* \right) \leq u_i \leq \mathop {\max} \left( u_{k,i}, u_{k+1,i}^* \right), \;\; i = 1,...,n_u  \},
\end{equation}

\noindent i.e. the box formed by the two points. Since we cannot leave this subspace between iterations $k$ and $k+1$, it follows that the Lipschitz constants may be found by the same means as in (\ref{eq:parlip}) but with the constraint $u_1,u_2 \in \mathcal{I}_k$, which may tighten the gap between $\underline \kappa_{j2}$ and $\overline \kappa_{j2}$ significantly.

We do note that numerical problems like (\ref{eq:parlip}) may not always be easily solved, however, and require reliable global optimization.

\subsubsection{Exploiting Partial Local Concavity}

A significant relaxation to the Lipschitz bound may be made if a plant constraint is assumed to be concave in some of the variables on the local input space $\mathcal{I}_k$.

\begin{thm}{\bf (Concavity-Based Relaxation of the Linear Upper Bound)}

Let $g_{p,j}$ be strictly concave\footnote{To stay in line with the previous theory of strict upper bounds generated by the Lipschitz constants, we will, for simplicity, assume \emph{strict} concavity.} in ${\bf v} = \{u_1,...,u_c\}, \; \forall {\bf u} \in \mathcal{I}_k$, with the remaining variables denoted by ${\bf z} = \{ u_{c+1},...,u_{n_u} \}$. Then, the following relaxation of (\ref{eq:linbound}) is valid:

\begin{equation}\label{eq:hybridbound}
\begin{array}{l}
g_{p,j} ( {\bf{u}}_{k+1}) - g_{p,j} ({\bf{u}}_k) < \vspace{1mm}\\
\hspace{10mm}  {\displaystyle \sum\limits_{i = 1}^{c} {\frac{\partial g_{p,j}}{\partial u_i} \Big |_{{\bf u}_k} ( u_{k+1,i} - u_{k,i} )}} + {\displaystyle \sum\limits_{i = c+1}^{n_u} {\kappa_{ji} |  u_{k+1,i} - u_{k,i} | }} 
\end{array}.
\end{equation}

\end{thm}

\begin{pf}

We start by decomposing the evolution of $g_{p,j}$ into the concave and non-concave portions:

\begin{equation}\label{eq:evodecomp}
\begin{array}{l}
g_{p,j} ( {\bf{u}}_{k+1}) - g_{p,j} ({\bf{u}}_k) = \vspace{1mm} \\
\left[ g_{p,j} ( {\bf{v}}_{k+1},{\bf{z}}_{k}) - g_{p,j} ( {\bf{v}}_{k},{\bf{z}}_{k}) \right] + \left[ g_{p,j} ( {\bf{v}}_{k+1}, {\bf{z}}_{k+1}) - g_{p,j} ( {\bf{v}}_{k+1},{\bf{z}}_{k}) \right] 
\end{array}.
\end{equation}

From strict concavity and from the partial version of (\ref{eq:linbound}), we have:

\begin{equation}\label{eq:sepbound}
\begin{array}{l}
g_{p,j} ( {\bf{v}}_{k+1},{\bf{z}}_{k}) - g_{p,j} ( {\bf{v}}_{k},{\bf{z}}_{k}) \leq {\displaystyle \sum\limits_{i = 1}^{c} {\frac{\partial g_{p,j}}{\partial u_i} \Big |_{{\bf u}_k} ( u_{k+1,i} - u_{k,i} )}} \\
g_{p,j} ( {\bf{v}}_{k+1}, {\bf{z}}_{k+1}) - g_{p,j} ( {\bf{v}}_{k+1},{\bf{z}}_{k}) \leq {\displaystyle \sum\limits_{i = c+1}^{n_u} {\kappa_{ji} |  u_{k+1,i} - u_{k,i} | }}
\end{array},
\end{equation}

\noindent with strict inequality for the top and bottom cases provided that ${\bf v}_{k+1} \neq {\bf v}_k$ and ${\bf z}_{k+1} \neq {\bf z}_k$, respectively. Since only ${\bf u}_{k+1} \neq {\bf u}_k$ is being considered (the linear upper bound not being of interest otherwise), it follows that at least one of the inequalities is strict, and that the two may then be added to yield (\ref{eq:hybridbound}). \qed

\end{pf}

Substituting in the input filter law and following a similar derivation to that of Theorem 2 in the companion work, it is simple to show that we may enforce $g_{p,j} ( {\bf{u}}_{k+1}) < 0$ via the following condition on $K_k$:

\begin{equation}\label{eq:liplimithybrid}
K_k \leq \frac{\displaystyle -g_{p,j} ({\bf u}_{k} ) }{{\displaystyle \sum\limits_{i = 1}^{c} {\frac{\partial g_{p,j}}{\partial u_i} \Big |_{{\bf u}_k} (  u^*_{k+1,i} - u_{k,i} )}} + \displaystyle \sum\limits_{i = c+1}^{n_u} {\kappa_{ji} |  u^*_{k+1,i} - u_{k,i} | }},
\end{equation}

\noindent provided that the denominator is positive (otherwise, the feasibility-guaranteeing $K_k$ has no upper limit and may be set to 1). It is clear that every additional input in which $g_{p,j}$ is concave will lead to an increase in the upper bound on $K_k$ by virtue of:

\begin{equation}\label{eq:lipvsccv}
\frac{\partial g_{p,j}}{\partial u_i} \Big |_{{\bf u}_k} ( u^*_{k+1,i} - u_{k,i} ) < \kappa_{ji} | u^*_{k+1,i} - u_{k,i} | , \;\;\; i = 1,...,c.
\end{equation}

The robust version of (\ref{eq:liplimithybrid}), which takes into account the uncertainty in the derivatives and the constraint measurements (as outlined in the previous sections) may be obtained by making the following substitutions:

\begin{equation}\label{eq:liplimithybridrob}
\begin{array}{l}
{\displaystyle \sum\limits_{i = 1}^{c} {\frac{\partial g_{p,j}}{\partial u_i} \Big |_{{\bf u}_k} (  u^*_{k+1,i} - u_{k,i} )}} \rightarrow \\
\displaystyle \sum\limits_{i = 1}^{c} {{\rm{max}} \left({\frac{\partial \overline g_{p,j}}{\partial u_i} \Big | _{{\bf{u}}_k} ( u^*_{k+1,i} - u_{k,i} ) }, {\frac{\partial \underline g_{p,j}}{\partial u_i} \Big | _{{\bf{u}}_k} ( u^*_{k+1,i} - u_{k,i} ) }\right)} \\
g_{p,j}({\bf u}_{k} ) \rightarrow \overline g_{p,j} ({\bf u}_{k} )
\end{array},
\end{equation}

\noindent where the bounded derivatives are the original robust bounds corresponding to the fully robust case, and not those obtained after reducing the bounds as described in Section 2.

While assuming plant data to have special structures like concavity is not new conceptually \citep{Ubhaya:09}, it is not a standard assumption and virtually absent in RTO, with the recent work in \cite{Bunin:12} appearing to be the first attempt to exploit such special structures in the RTO context. This does not, however, preclude the existence of cases where such assumptions would be justified. Taking again the example in (\ref{eq:paruncertainty}), we see that such an assumption could be made if we knew that the parameter $\theta$ were negative (since the constraint would then be globally concave in $u_2$). Alternatively, one could use available plant data to test if a concave model might be an appropriate fit for the data in a given neighborhood -- such tests are tractable and can be done online \citep[Section 6.5]{Boyd:08}\citep{Bunin2013a}. While potentially more risky, the latter presents a data-driven way to help validate, or even ``discover'', such an assumption.

\vspace{2mm}
\noindent{\it{Example -- Assumption of Concavity for Problem A}} 
\vspace{2mm}

\begin{figure}
\begin{center}
\includegraphics[width=8cm]{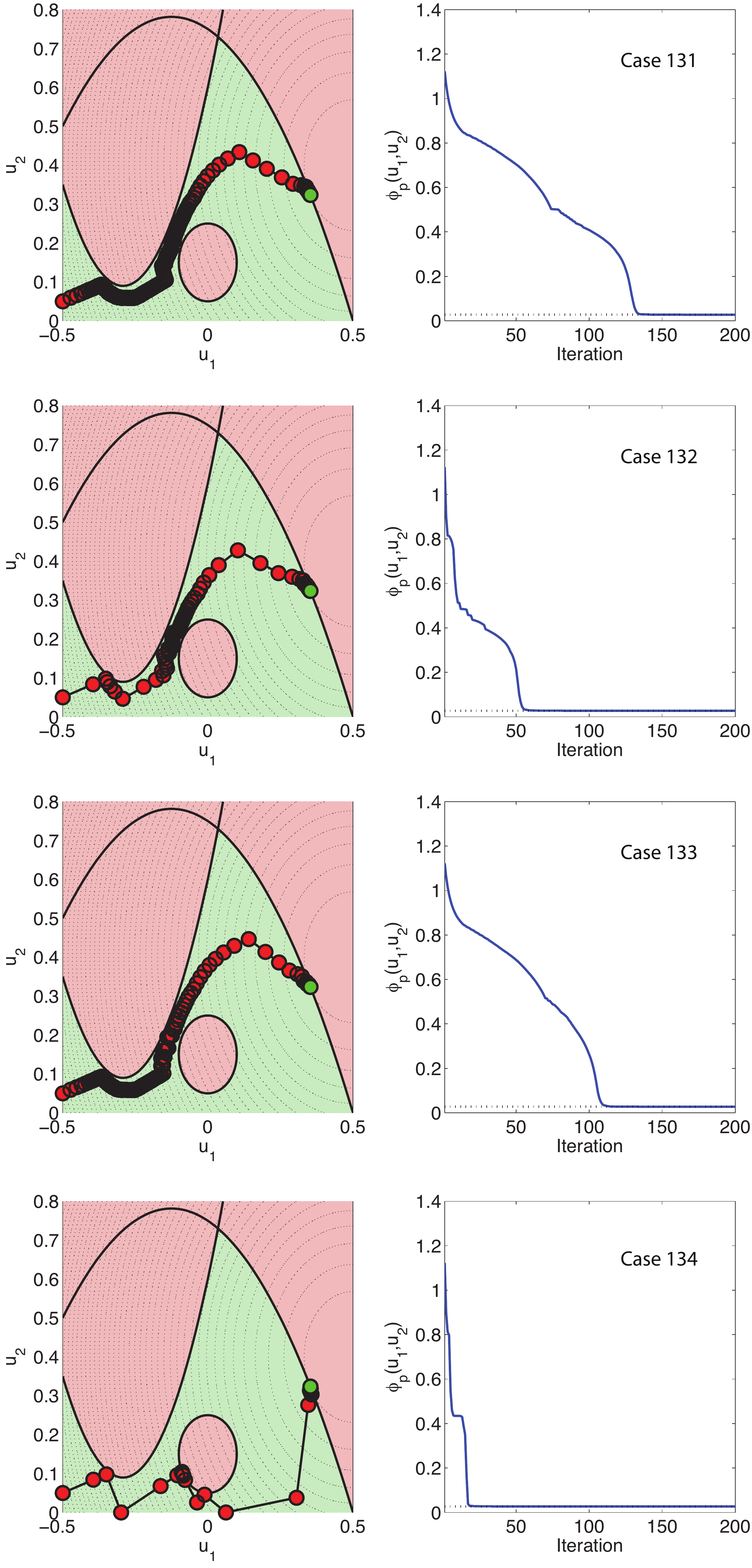}    
\caption{Improvement in RTO performance that results from assuming partial concavity in the constraint functions.}
\label{fig:P2ex10}
\end{center}
\end{figure}

We again consider Problem A, for which we know that the first and third constraints are globally concave, and examine the performance benefits brought by assuming this structure formally. Namely, we compare the performance of the ideal-target algorithm over 200 iterations when: (a) no concavity is assumed and the standard Lipschitz constants are used throughout, (b) concavity with respect to $u_1$ is assumed for the two concave constraints, (c) concavity with respect to $u_2$ is assumed for the two concave constraints, and (d) concavity with respect to both inputs is assumed for both constraints. The performances are given in Table \ref{tab:concave} and in Figure \ref{fig:P2ex10}, which demonstrate how such an assumption on the structure of the constraint function may lead to dramatic improvements in performance.

\begin{table}
\caption{Performance of the ideal-target scheme with different assumptions regarding constraint concavity for Problem A.}
\begin{center}
\footnotesize{
\begin{tabular}{|p{.6cm}|p{5cm}|p{1cm}|}
\hline
Case & Assumption & $L$   \\ \hline
131 & None & 73.54 \\ \hline
132 & $g_{p,1}$ and $g_{p,3}$ concave w.r.t. $u_1$ & 22.16 \\ \hline 
133 & $g_{p,1}$ and $g_{p,3}$ concave w.r.t. $u_2$ & 63.46 \\ \hline 
134 & $g_{p,1}$ and $g_{p,3}$ concave w.r.t. $u_1, u_2$ & 8.17 \\ \hline
\end{tabular}
}
\end{center}
\label{tab:concave}
\end{table}

\subsubsection{Using Previous Measurements}

As illustrated in the companion work, the Lipschitz bound around a strictly feasible input point essentially generates a polytope that is guaranteed to contain only strictly feasible points with regard to $g_{p,j}$:

\begin{equation}\label{eq:lippoly}
\mathcal{L}_{k,j} = \left\{ {\bf u} : g_{p,j} ({\bf{u}}_k) + \displaystyle \sum\limits_{i = 1}^{n_u} {\kappa_{ji} | u_{i} - u_{k,i} | } \leq 0 \right\},
\end{equation}

\noindent which can, of course, be refined further with the ideas of Sections 4.2.1 and 4.2.4, and made robust with substitutions like those in (\ref{eq:liplimithybridrob}).

Since every single data point available from prior operation will generate its own $\mathcal{L}$, we may relieve some of the locality of the Lipschitz bound by observing that the union of these polytopes $\mathcal{L}_{U,j}$:

\begin{equation}\label{eq:polyint}
\mathcal{L}_{U,j} = \bigcup_{i = 0}^{k} \mathcal{L}_{i,j},
\end{equation}

\noindent must contain only strictly feasible points as well.

This may become useful whenever the RTO algorithm revisits a certain region of the input space. While an ideal algorithm may be imagined to take a fairly straight trajectory towards the optimum (and thereby not visit the same region twice), there are some practical cases when such behavior could occur:

\begin{itemize}
\item The algorithm does not perform well due to significant uncertainty in the gradients (which could lead to haphazard steps).
\item The algorithm performs well but converges in a zig-zag manner.
\item The algorithm converges, but the RTO problem changes (e.g. due to a change in user/market demand, which may change the cost function), and the algorithm must converge to a new optimum.
\item Experimental data from the process is available prior to the launching of the RTO algorithm.
\end{itemize}

We may then relax the requirement on $K_k$ by finding it in the following manner:

\begin{equation}\label{eq:linesearchK}
\begin{array}{l}
\mathop {{\rm{maximize}}}\limits_{K_k, {\bf u}_{k+1}}\hspace{3mm}K_k \\
{\rm{subject}}\hspace{1mm}{\rm{to}}\hspace{3mm}{\bf{u}}_{k+1}  = {\bf{u}}_k  + K_k \left( \bar {\bf{u}}^*_{k + 1} - {\bf{u}}_k \right) \\
\hspace{18mm}{\bf u}_{k+1} \in \displaystyle \bigcap_{j=1}^{n_g} \mathcal{L}_{U,j} \\
\hspace{18mm}K_k \leq 1
\end{array},
\end{equation}

\noindent that is, we allow potentially bigger steps if any values of $K_k$ yield inputs belonging to feasible polytopes from past iterations. While the feasible space of (\ref{eq:linesearchK}) may be both disjoint and/or nonconvex, this problem is easy to solve as it only involves a line search in $K_k$.

\begin{figure}
\begin{center}
\includegraphics[width=8cm]{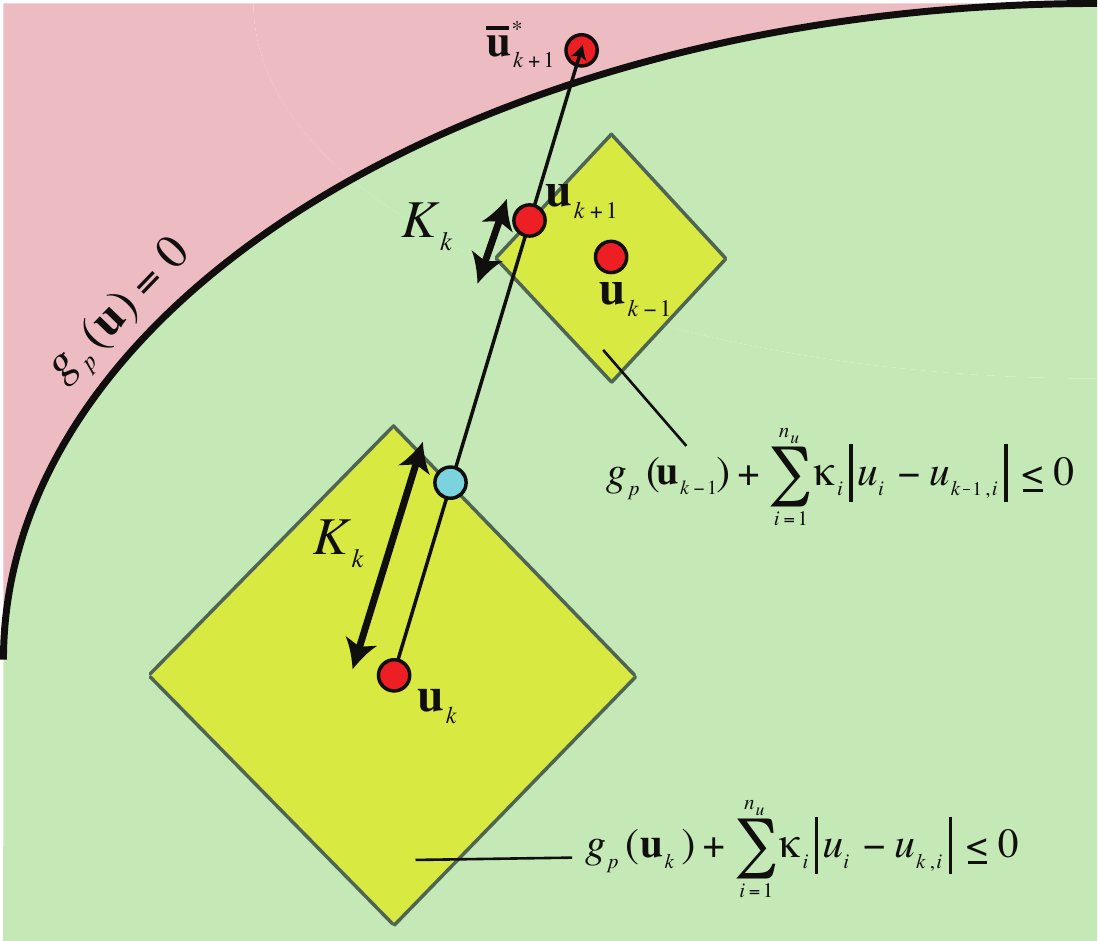}    
\caption{Past data may be used to allow for longer input steps. Here, the optimal direction goes through the feasible polytope around ${\bf u}_{k-1}$, which allows the algorithm to take a larger step than it normally would (the latter given by the blue dot). The two-sided arrows show all portions of the line search (all values of $K_k$) that are guaranteed to yield feasible iterates.}
\label{fig:lippast}
\end{center}
\end{figure}

A geometric illustration of this simple concept is given in Figure \ref{fig:lippast}.

\vspace{2mm}
\noindent{\it{Example -- Problem B with a Changing Cost Function}} 
\vspace{2mm}

To illustrate the above ideas, we consider Problem B and run it for 100 iterations. For the first 49 iterations, the standard definition of the RTO problem is used. From Iteration 50 onwards, we change the cost function of the problem to:

\begin{equation}\label{eq:costchange}
\phi_{p}({\bf{u}}) = (u_1+0.25)^2 + (u_2-0.6)^2,\;\;{\rm when}\; k \geq 50,
\end{equation}

\noindent which then forces the RTO algorithm to go to the new optimum.

\begin{figure}
\begin{center}
\includegraphics[width=7cm]{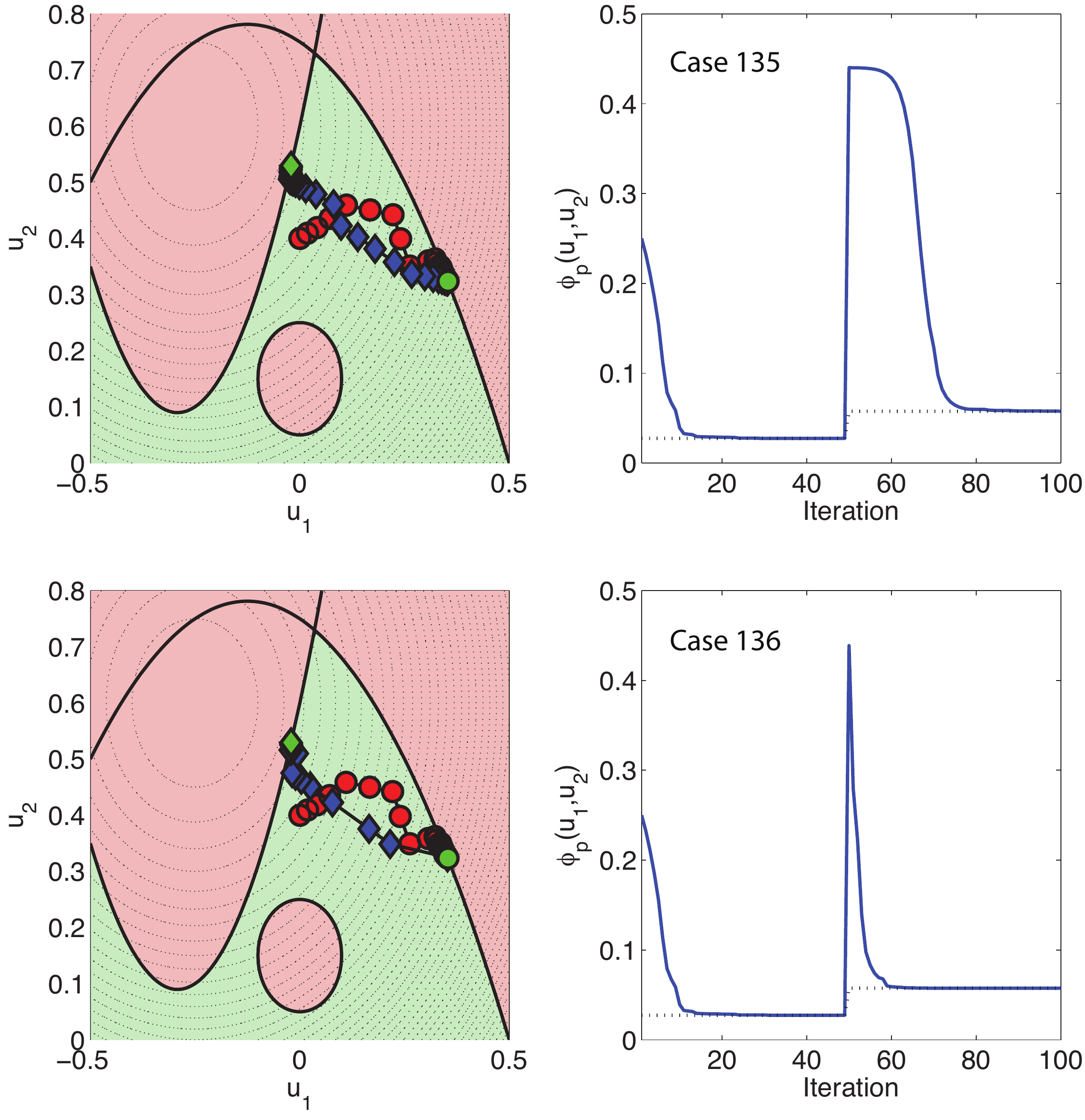}    
\caption{The use of past data points to relax the Lipschitz bound can yield a significant improvement in performance for a problem where the cost function changes during operation. Only the contours for the second cost function are given.}
\label{fig:P2ex11}
\end{center}
\end{figure}

We test this for two cases -- one where we do not attempt to relax the Lipschitz bound by using past measurements and one where we do. The results are presented in Table \ref{tab:pastmeas} and in Figure \ref{fig:P2ex11}, with the first converging sequence given by the usual red dots and the second given by blue diamonds (with the green diamond denoting the plant optimum for the second cost function). We see that, while the first sequence is the same for both cases (the algorithm does not visit the same region twice), the second converges to the second optimum significantly faster when past measurement knowledge is incorporated to relax the bound on $K_k$. Clearly, this is because the algorithm is forced to cross through regions where it has been before.

\begin{table}
\caption{Performance of the ideal-target scheme when past measurements are exploited to relax the Lipschitz bound.}
\begin{center}
\footnotesize{
\begin{tabular}{|p{.6cm}|p{4cm}|p{1cm}|}
\hline
Case & Past Measurements Used? & $L$   \\ \hline
135 & No & 7.81 \\ \hline
136 & Yes & 2.12 \\ \hline 
\end{tabular}
}
\end{center}
\label{tab:pastmeas}
\end{table}

\subsection{Quadratic Upper Bound on the Cost}

A lot of the concepts that apply to the Lipschitz constants extend quite naturally to the quadratic upper bound, which we note is essentially derived by using the Lipschitz property of the \emph{second-order derivatives} of a function (see Lemma 2 in the companion work). However, as experimental second-order derivatives are significantly more difficult to estimate, or even derive, than the first-order ones, the idea of a ``physical significance'' is expected to be absent for $\overline {\bf Q}_{\phi}$, and the techniques used to relax or refine $\overline {\bf Q}_{\phi}$ are expected to be less reliable. Because of these challenges, we start by discussing when a good choice of $\overline {\bf Q}_{\phi}$ is crucial, as cases where this is not so may not cause problems for the user. We then go through the different refinement techniques, skipping the derivations and stating only the results, as many are identical in form to what has already been given for the first-order Lipschitz cases. We finish by suggesting a general-purpose adaptive choice of $\overline {\bf Q}_{\phi}$ that would guarantee convergence as $k \rightarrow \infty$.

\subsubsection{When Does the Choice of $\overline {\bf Q}_{\phi}$ Matter?} 

Recall that a $\overline {\bf Q}_{\phi}$ with values that are too small can cause the RTO algorithm to take steps that do not guarantee a decrease in the cost at each iteration, which in turn loses the conceptual guarantee of finite-time convergence to a KKT point. On the other hand, defining $\overline {\bf Q}_{\phi}$ with values that are too large can make the algorithm move very slowly (due to Condition (\ref{eq:suffcond4})), which causes performance loss due to slow convergence. Here, we are mainly interested in the former as it illustrates, at least from the theoretical point of view, a more detrimental case where the algorithm is not able to converge to an optimum even as $k \rightarrow \infty$. As such, we can rephrase the question above as ``when does selecting a $\overline {\bf Q}_{\phi}$ that does not satisfy (\ref{eq:Qbound}) have a crucial effect on performance?''. 

We give a quantitative, though abstract, answer to this question by stating that whenever

\begin{equation}\label{eq:Qqual}
\mathop{\min}\limits_{j = 1,...,n_g} \left[ \frac{\displaystyle - g_{p,j} ({\bf u}_{k} ) }{\displaystyle \sum\limits_{i = 1}^{n_u} {\kappa_{ji} | \bar u^*_{k+1,i} - u_{k,i} | }} \right] \leq -2 \displaystyle \frac{\nabla \phi_p({\bf{u}}_k)^T ( \bar {\bf{u}}^*_{k+1} - {\bf{u}}_{k})}{( \bar {\bf{u}}^*_{k+1} - {\bf{u}}_{k})^T \overline {\bf{Q}}_\phi (\bar {\bf{u}}^*_{k+1} - {\bf{u}}_{k})},
\end{equation}

\noindent the choice of $\overline {\bf Q}_{\phi}$ becomes irrelevant because it does not define the upper bound on the input filter gain (the latter being defined by the feasibility condition). We may also consider the robust version of the above:

\begin{equation}\label{eq:Qqualrob}
\begin{array}{l}
\mathop{\min}\limits_{j = 1,...,n_g} \left[ \frac{\displaystyle - \bar g_{p,j} ({\bf u}_{k} ) }{\displaystyle \sum\limits_{i = 1}^{n_u} {\kappa_{ji} | \bar u^*_{k+1,i} - u_{k,i} | }} \right] \leq \vspace{2mm} \\
-2\frac{\displaystyle \sum\limits_{i = 1}^{n_u} {{\rm{max}} \left({\frac{\partial \overline \phi_p}{\partial u_i} \Big | _{{\bf{u}}_k} (\bar u^*_{k+1,i} - u_{k,i} ) }, {\frac{\partial \underline \phi_p}{\partial u_i} \Big | _{{\bf{u}}_k} (\bar u^*_{k+1,i} - u_{k,i} ) }\right)}}{\displaystyle (\bar {\bf{u}}^*_{k+1} - {\bf{u}}_{k})^T \overline {\bf{Q}}_\phi (\bar {\bf{u}}^*_{k+1} - {\bf{u}}_{k})}
\end{array}.
\end{equation}

From this, we may make a number of qualitative observations regarding the cases where (\ref{eq:Qqualrob}) would be satisfied:

\begin{itemize}
\item When there are many constraints (large $n_g$).
\item When the uncertainty in the constraint measurements/estimates is large and/or the operating point is close to some constraints (small $-\overline g_{p,j}({\bf u}_{k} )$).
\item When the Lipschitz constants are very conservative (large) and few attempts are made to relax them (using the methods of the previous subsection, for example).
\item When there is little uncertainty in the gradients of the cost function (as uncertainty would automatically increase the numerator of the right-hand side, thereby lowering the value of the bound).
\end{itemize}

The overall statement, and not a surprising one, is that a false choice of $\overline {\bf Q}_{\phi}$ is not as debilitating when the RTO problem is highly constrained, as this contributes to lowering the bound on the filter gain to enforce feasibility guarantees. This has indeed been observed by us in the numerical trials for Problem A, where the upper bound on the filter gain given by Condition (\ref{eq:suffcond4}) was never observed to be active, and where the $K_k$ was always bounded by Condition (\ref{eq:suffcond3}). By contrast, an unconstrained problem would rely entirely on the proper choice of $\overline {\bf Q}_{\phi}$ so as not to lose convergence guarantees (depending on the choice of algorithm).

Our qualitative advice is therefore to take more care with the choice of $\overline {\bf Q}_{\phi}$ when there are few constraints in the problem and when there is a reasonable chance of the optimum being unconstrained. That being said, we now present several methods to make this choice easier or to refine it, many of which are analogous to the choice and refinement of the Lipschitz bounds.

\subsubsection{Choosing and Refining the Quadratic Upper Bound}

Referring the reader to Lemma 2 and its derivation in the companion work, we start by proposing the following relaxation of the second-order Lipschitz constants $M$:

\begin{equation}\label{eq:Mrel}
\underline M_{ij} < \frac{\partial ^2 \phi_p}{\partial u_i \partial u_j} \Big |_{\bf u} < \overline M_{ij}, \;\; \forall {\bf u} \in \mathcal{I},
\end{equation}

\noindent i.e. the second-order analogue of (\ref{eq:liptwo}). Without going through the derivation, we note that all of the previously stated concepts, like those of (\ref{eq:der2}), may be reused here to yield the following relaxation:

\begin{equation}\label{eq:Qrel}
\begin{array}{l}
({\bf u}_{k+1} - {\bf u}_{k})^T \overline {\bf Q}_{\phi} ({\bf u}_{k+1} - {\bf u}_{k}) \geq \\
\hspace{20mm} \displaystyle \sum_{i=1}^{n_u} \sum_{j=1}^{n_u} \mathop {\max} \Big ( \underline M_{ij} (u_{k+1,i} - u_{k,i} )(u_{k+1,j} - u_{k,j} ), \\
\hspace{40mm}\overline M_{ij} (u_{k+1,i} - u_{k,i} )(u_{k+1,j} - u_{k,j} )  \Big )
\end{array},
\end{equation}

\noindent with the right-hand side substituted for the left-hand side in (\ref{eq:suffcond4}) to relax the filter gain bound accordingly (provided the right-hand side is positive).

We may then take advantage of this characterization by noting that the absence of nonlinearity between certain variables and the cost allows $\underline M_{ii}, \overline M_{ii} \approx 0$ (i.e. negative and positive numbers that are very close to 0, so as to enforce the strictness of (\ref{eq:Mrel})), and that the absence of interaction allows $\underline M_{ij} = \underline M_{ji} \approx 0$ and $\overline M_{ij} = \overline M_{ji} \approx 0$. This is analogous to exploiting the sparsity in the Jacobian when defining the first-order Lipschitz constants.

We may also exploit any uncertainty that is known to be exclusively parametric. If, for example, we know that the cost can be properly described by:

\begin{equation}\label{eq:paruncertaintycost}
\phi ({\bf u},{\boldsymbol \theta}) = \theta_1 u^2_1 + \theta_2 u_2^2 + \theta_1 \theta_2 u_1 u_2,
\end{equation}

\noindent with $\theta_1 \in (\underline \theta_1, \overline \theta_1)$ and $\theta_2 \in (\underline \theta_2, \overline \theta_2)$, then we can analytically calculate the second derivatives:

\begin{equation}\label{eq:secondder}
\frac{\partial ^2 \phi}{\partial u_1^2} = 2\theta_1,\;\; \frac{\partial ^2 \phi}{\partial u_1 \partial u_2} = \frac{\partial ^2 \phi}{\partial u_2 \partial u_1} = \theta_1 \theta_2,\;\; \frac{\partial ^2 \phi}{\partial u_2^2} = 2\theta_2,
\end{equation}

\noindent and then define the $M$ constants as:

\begin{equation}\label{eq:secondderM}
\begin{array}{l}
\underline M_{11} = 2 \underline \theta_1, \overline M_{11} = 2 \overline \theta_1 \\
\underline M_{12} = \underline M_{21} = \mathop {\min} \{ \theta_1 \theta_2 : \underline \theta_1 \leq \theta_1 \leq \overline \theta_1, \underline \theta_2 \leq \theta_2 \leq \overline \theta_2 \} \\
\overline M_{12} = \overline M_{21} = \mathop {\max} \{ \theta_1 \theta_2 : \underline \theta_1 \leq \theta_1 \leq \overline \theta_1, \underline \theta_2 \leq \theta_2 \leq \overline \theta_2 \} \\
\underline M_{22} = 2 \underline \theta_2, \overline M_{22} = 2 \overline \theta_2
\end{array},
\end{equation}

\noindent which would then be used in (\ref{eq:Qrel}).

If the cost is known to be concave in a certain variable, then it follows that $\overline M_{ii} \leq 0$. Likewise, $\underline M_{ii} \geq 0$ would result if the relationship were convex.

We can also exploit past measurements. Introducing $\mathcal{F}_{i}$ and defining it as:

\begin{equation}\label{eq:Fpoly}
\begin{array}{l}
\mathcal{F}_{i} = \Big \{ {\bf u} : \phi_p ({\bf{u}}_{i}) + \nabla \phi_p ({\bf{u}}_{i})^T ({\bf{u}} - {\bf{u}}_{i}) + \\
\hspace{15mm}\displaystyle \frac{1}{2}({\bf{u}} - {\bf{u}}_{i})^T \overline {\bf{Q}}_\phi ({\bf{u}} - {\bf{u}}_{i}) < \phi_p ({\bf{u}}_{k}) \Big \},
\end{array}
\end{equation}

\noindent we may generate a union of such quadratic sets:

\begin{equation}\label{eq:Funion}
\mathcal{F}_U = \bigcup_{i=0}^{k}\mathcal{F}_i,
\end{equation}

\noindent all of which are guaranteed to contain points with cost function values strictly inferior to the value at the current iterate. We may then attempt to find a larger value of $K_k$ by employing a line search like that of (\ref{eq:linesearchK}).

\subsubsection{A General Algorithm to Determine $\overline {\bf Q}_\phi$ Online}

In the interest of simplicity, we propose the following algorithm to adapt $\overline {\bf Q}_\phi$ iteratively when $\overline {\bf Q}_\phi \succ 0$:

\begin{enumerate}
\item Set $n := 0$, $\overline {\bf Q}_\phi := {\bf I}$, $k :=0$.
\item RTO Iteration: $k := k + 1$.
\item If $\phi_p ({\bf u}_k) \geq \mathop{\min}\limits_{i = n,...,k-1} \phi_p ({\bf u}_i)$ with sufficiently high confidence, set $\overline {\bf Q}_\phi := 2 \overline {\bf Q}_\phi$ and $n := k$.
\item Return to Step 2.
\end{enumerate}

It is easy to show that such an algorithm preserves the guarantee of finite-time convergence since consistent failure to decrease $\phi_p$ monotonically would force $\overline {\bf Q}_\phi \rightarrow 2^\infty {\bf I}$, which is clearly a sufficiently high upper bound. A price to pay for such a simple implementation is that of potentially augmenting the cost in the beginning. If this is a major concern, then using something greater than ${\bf I}$ to initialize the algorithm could also be proposed.

\section{Improving Convergence Speed}

In this final part, we consider two important scenarios that may improve the convergence speed of the SCFO-supplemented RTO algorithm. We consider:

\begin{itemize}
\item allowing constraint violations during convergence,
\item assuming that some of the constraints and/or the cost are known exactly,
\end{itemize}

\noindent and show how these scenarios may be incorporated into the SCFO. For simplicity, we work with the nominal formulations, but note that all of the results that follow may be made robust with substitutions like (\ref{eq:liplimithybridrob}), and that the relaxed formulations of Section 4 may also be included with ease.

\subsection{Allowing Constraint Violations}

Up to now, we have considered all of the constraints in (\ref{eq:realopt}) as being hard, i.e. as having to be satisfied at \emph{every} RTO iteration. In practice, however, one is often confronted with soft constraints as well, which may be violated at certain iterations and only need to be satisfied upon convergence. We may incorporate the idea of a soft constraint $g_{p,j}$ into our framework by:

\begin{enumerate}
\item Specifying the maximal allowable violation for this constraint and denoting it by a strictly positive value $d_j$.
\item Specifying a maximal violation \emph{integral}, $d_{T,j}$, or the total sum of violations that we are willing to tolerate over \emph{all} RTO iterations, which leads to the constraint:

\begin{equation}\label{eq:violint}
\displaystyle \sum\limits_{i = 0}^{\infty} {\rm max} (0,g_{p,j}({\bf u}_i)) \leq d_{T,j}.
\end{equation}

\end{enumerate}

As such, we would allow our RTO algorithm to violate, in the worst case, $g_{p,j}$ by $d_j$, but would prefer that it not do so indefinitely, and that it return to the feasible region before the limit (\ref{eq:violint}) is passed.

The first of these requirements may be incorporated quite readily into the filter gain bound (we skip the derivation as it is akin to that of Theorem 2 in \cite{Bunin:12b}):

\begin{equation}\label{eq:gainupperslack}
K_{k} \le \mathop{\min}\limits_{j = 1,...,n_g} \left[ \frac{{- g_{p,j}({\bf{u}}_{k } ) + d_j}}{\displaystyle \sum\limits_{i = 1}^{n_u} {\kappa_{ji} | u^*_{k+1,i} - u_{k,i} | }} \right],
\end{equation}

\noindent which has the performance-enhancing property of allowing larger steps and generally faster convergence.

To enforce (\ref{eq:violint}), we propose to gradually decrease $d_j$ via the following slack reduction law for every iteration $k$ where $ g_{p,j}({\bf{u}}_k) \geq 0$:

\begin{equation}\label{eq:violupdate}
d_j := \beta_j d_j,
\end{equation}

\noindent with $\beta_j \in [0,1)$ the slack reduction law constant.

We now provide a choice of $\beta_j$ that is sufficient for guaranteeing (\ref{eq:violint}).

\begin{thm}{\bf (Upper Bound on Slack Reduction Law Constant)}

Denote by $I_j$ the set of iteration indices, ordered from smallest to largest, where $g_{p,j}({\bf{u}}_k) \geq 0$, i.e. $I_j = \left\{ i : g_{p,j}({\bf{u}}_i) \geq 0 \right\}$. If (\ref{eq:violupdate}) is applied for every $k \in I_j$, and the following limit on $\beta_j$ is met:

\begin{equation}\label{eq:betalimit}
\beta_j \leq \frac{d_{T,j}-d_{j,0}}{d_{T,j}},
\end{equation}

\noindent where $d_{j,0}$ is the initial $d_j$ value, then (\ref{eq:violint}) is guaranteed.

\end{thm}

\begin{pf}

We may rewrite the violation integral as:

\begin{equation}\label{eq:violintre}
\displaystyle \sum\limits_{i = 0}^{\infty} {\rm max} (0,g_{p,j}({\bf u}_i)) = \displaystyle \sum\limits_{i = 0}^{\| I_j \|_0 - 1} g_{p,j}({\bf u}_{I_j (i)}),
\end{equation}

\noindent with $I_j(0)$ denoting the first element of $I_j$.

From the bound (\ref{eq:gainupperslack}) and the update (\ref{eq:violupdate}), we have the guarantee that:

\begin{equation}\label{eq:slackg}
g_{p,j} ({\bf u}_{I_j (i)}) < \beta_j^i d_{j,0},
\end{equation}

\noindent which, in summation, yields:

\begin{equation}\label{eq:violintsum}
\displaystyle \sum\limits_{i = 0}^{\| I_j \|_0 - 1} g_{p,j}({\bf u}_{I_j (i)}) < d_{j,0} \displaystyle \sum\limits_{i = 0}^{\| I_j \|_0 - 1} \beta_j^i = d_{j,0} \frac{1 - \beta_j^{\| I_j \|_0}}{1 - \beta_j} \leq \frac{d_{j,0}}{1 - \beta_j}.
\end{equation}

Forcing the right-hand side to be below the desired threshold, $d_{T,j}$, and rearranging yields (\ref{eq:betalimit}). \qed

\end{pf}

Practically, it may occur that the reduction in the slack is too fast and makes the iterate at $k$ infeasible with respect to the lowered slack at $k+1$. As a simple illustration of this, imagine that a given constraint is allowed to be violated by a value of 4, and that we would like for the total violation integral to remain below 6, which leads to $\beta_j \leq \frac{1}{3}$ if we employ (\ref{eq:betalimit}). Now, imagine that the first constraint violation takes us straight to a value of 3. Applying the slack reduction here would lower the slack from 4 to $\frac{4}{3}$ (or less), which would mean that the next iteration would need to guarantee a constraint value of less than $\frac{4}{3}$, despite the value already being at 3 at the current iteration. As the SCFO require always being on the feasible side of the constraint, we may handle this case with the following conditional:

\begin{itemize}
\item If $g_{p,j} ({\bf u}_k) \geq \beta_j d_j$, then find ${\bf u}_{\hat k}$, where $\hat k$ is a past iteration having the best observed cost and satisfying $g_{p,j} ({\bf u}_{\hat k}) < \beta_j d_j$. Using ${\bf u}_{\hat k}$ as a reference -- i.e. acting as if it were the current iterate -- carry out a standard RTO iteration with the new slack to determine ${\bf u}_{k+1}$.
\end{itemize}

\noindent This basically ensures that: if the algorithm is suddenly in violation of the constraint with the new slack, then the next iterate is simply obtained from the most optimal previous iterate where this constraint was met (note that such an iterate \emph{always} exists since the initial point, ${\bf u}_0$, is by assumption feasible for the 0-slack case). Such a scheme retains all of the original guarantees since continued constraint violations will eventually reduce the slack to 0 and, as the algorithm is structured to satisfy the slack constraint always, it will be forced to be feasible upon complete slack reduction as well. Furthermore, it is easy to show that convergence will be monotonic in the cost for all iterates that are feasible with respect to the original 0-slack constraints.

\vspace{2mm}
\noindent{\it{Examples with Soft Constraints}} 
\vspace{2mm}

\begin{figure}
\begin{center}
\includegraphics[width=8cm]{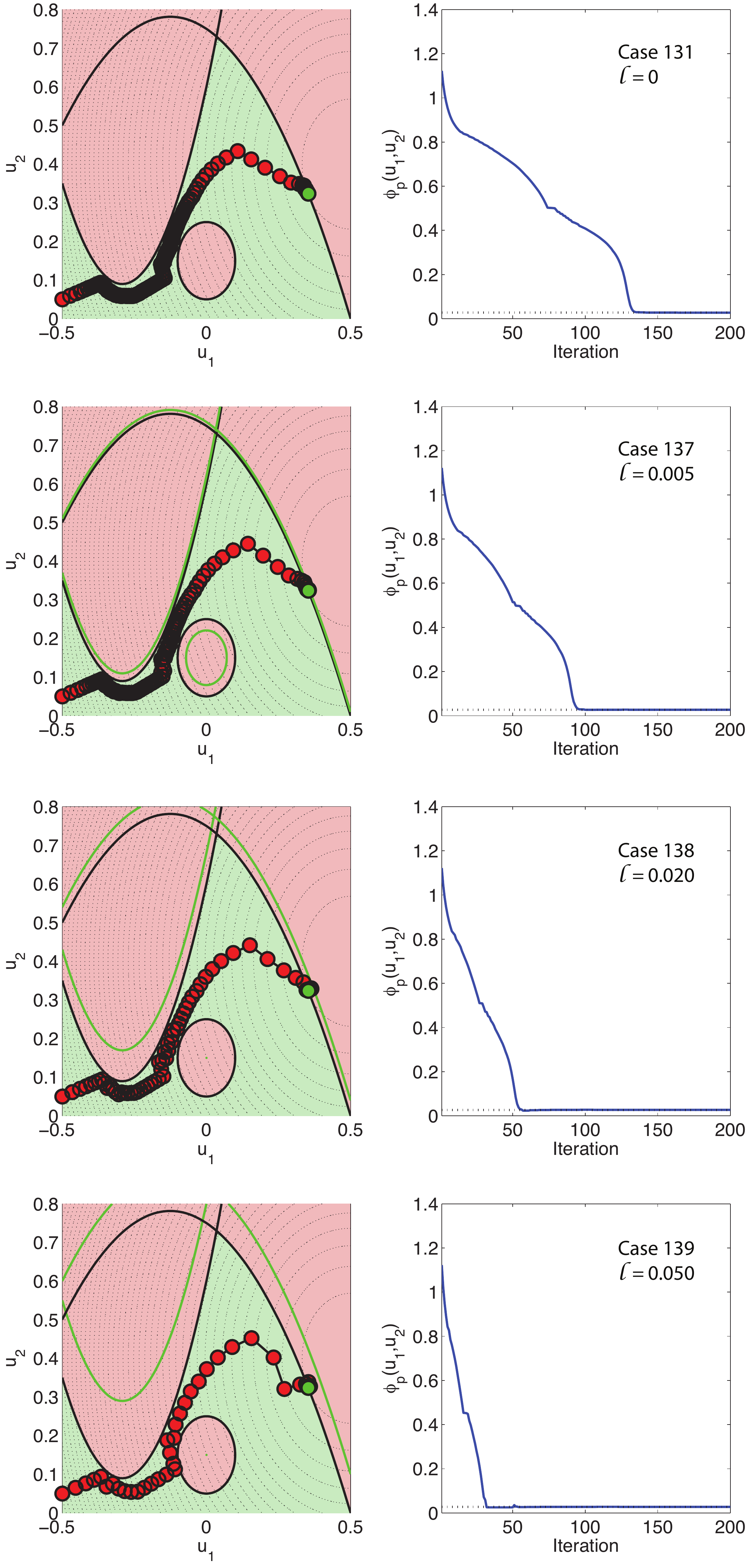}    
\caption{Performance of the ideal-target algorithm for Problem A when the constraints are made soft. Green lines denote the constraints with the added slack prior to slack reduction.}
\label{fig:P2ex12}
\end{center}
\end{figure}

\begin{figure}
\begin{center}
\includegraphics[width=8cm]{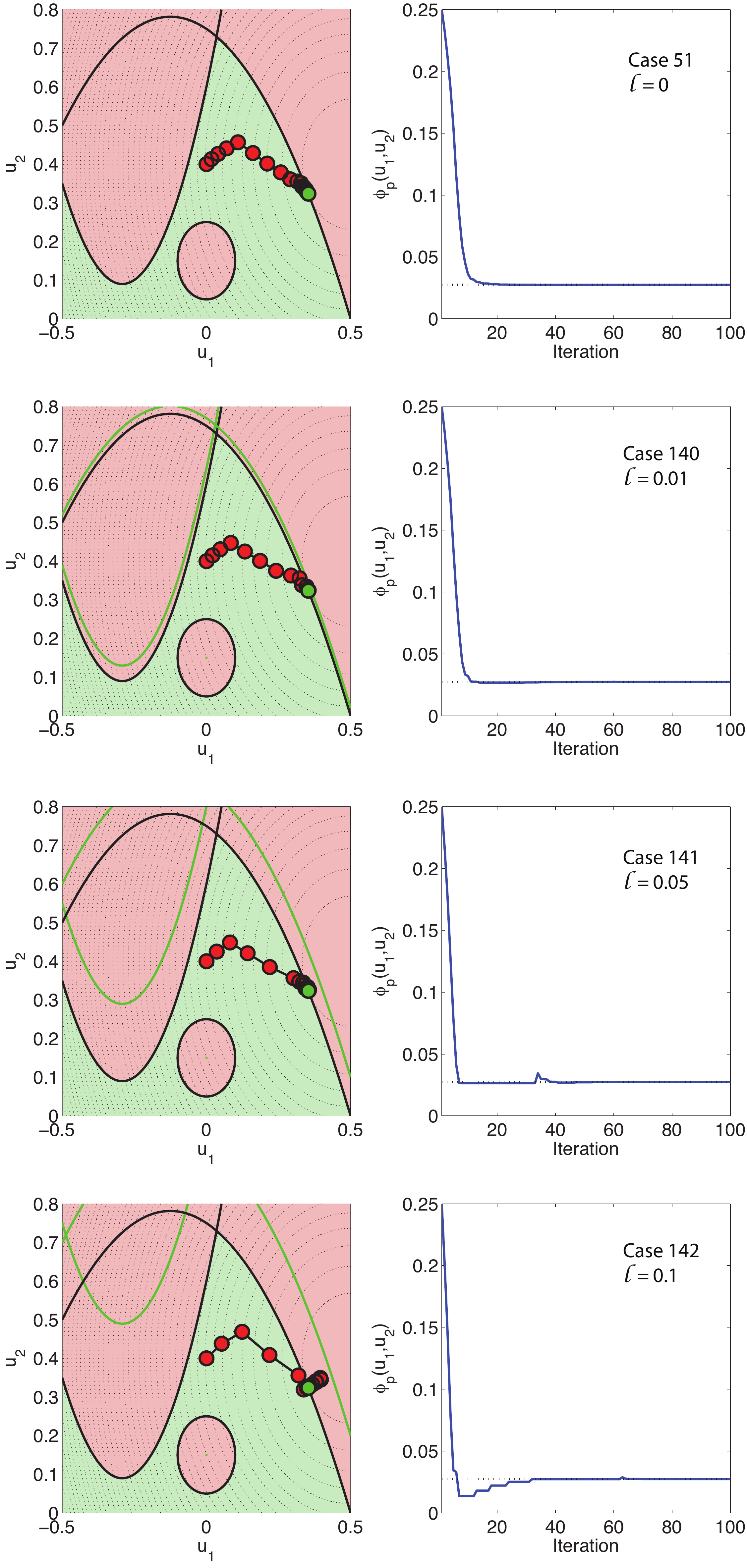}    
\caption{Performance of the ideal-target algorithm for Problem B when the constraints are made soft.}
\label{fig:P2ex13}
\end{center}
\end{figure}

We consider the same two problems as before, but this time soften the constraints proportionally to their ranges by specifying the maximal allowable violations and violation integrals as $d_{j,0} = l \overline \epsilon_j, d_{T,j} = 10 d_{j,0}$, where $l$ is a parameter that we can vary to augment or decrease the slack for testing purposes. As the results all follow the same trend (increasing $l$ increases performance), we only consider the ideal-target algorithm for each problem and report the results in Table \ref{tab:slackcon} and Figures \ref{fig:P2ex12} and \ref{fig:P2ex13}. The results are very encouraging as they demonstrate that not only may a dramatic improvement be achieved by softening the constraints, but that doing so often does not lead to serious violations -- as may be seen in Figure \ref{fig:P2ex12}, no violations at all were incurred in the concave constraints for Problem A even with significant slack. This is most likely due to the way in which the SCFO are implemented, as the projection innately attempts to avoid constraints as long as it does not need to approach them in order to lower the cost. Violations in the active constraint are, for this reason, inevitable, but we see that we are always able to reduce the slack to ultimately converge to the plant optimum.  

\begin{table}
\caption{Performance for different RTO problems when soft constraints are used.}
\begin{center}
\footnotesize{
\begin{tabular}{|p{.6cm}|p{1.2cm}|p{.6cm}|p{1cm}|}
\hline
Case & Problem & $l$ & $L$   \\ \hline
131 & A & 0 & 73.54 \\ \hline
137 & A & 0.005 & 50.78 \\ \hline 
138 & A & 0.020 & 27.93 \\ \hline 
139 & A & 0.050 & 15.46 \\ \hline
51 & B & 0 & 1.12 \\ \hline
140 & B & 0.010 & 0.98 \\ \hline 
141 & B & 0.050 & 0.71 \\ \hline 
142 & B & 0.100 & 0.40 \\ \hline 
\end{tabular}
}
\end{center}
\label{tab:slackcon}
\end{table}

\subsection{Incorporating Known Cost and Constraint Elements}

In many applications, it may happen that parts of Problem (\ref{eq:realopt}) are sufficiently well-modeled that we can allow the assumption of knowing those parts of the problem exactly. Alternatively, there may be problems where either the cost or the constraints are solely defined by the user, and as such are known analytically without uncertainty present. While all of these cases can still be solved by the general methods proposed in this work, we can gain a significant improvement in performance by incorporating this knowledge, as doing this not only avoids using upper bounds that are innately conservative, but may also remove constraints from the projection problem while retaining the desired guarantees (thereby diverting the target optimum less from that which is given by the RTO algorithm). We will break this discussion into three parts, based on the portions of the RTO problem that are known.

\subsubsection{Incorporating a Known Cost Function}

We consider the scenario where the cost is known but the constraints are not\footnote{We acknowledge knowledge of a function by removing the $p$ subscript.}:

\begin{equation}\label{eq:costknown}
\begin{array}{l}
\mathop {{\rm{minimize}}}\limits_{\bf{u}}\hspace{4mm}\phi ({\bf{u}}) \\
{\rm{subject}}\hspace{1mm}{\rm{to}}\hspace{3mm}{\bf{G}}_p({\bf{u}}) \preceq {\bf{0}} \\
\hspace{18mm}{\bf{u}}^L \preceq {\bf{u}} \preceq {\bf{u}}^U
\end{array}\hspace{3mm}.
\end{equation}

The initial step for solving this problem is to forego Condition (\ref{eq:suffcond1}) entirely, and to project only with respect to the uncertain constraints:

\begin{equation}\label{eq:feasprojcostknown}
\begin{array}{l}
\bar {\bf{u}}^*_{k+1} = {\rm{arg}} \mathop {{\rm{minimize}}}\limits_{{\bf{u}}}\hspace{4mm}\left\| {\bf{u}} - {\bf{u}}^*_{k+1} \right\|_2^2 \vspace{1mm}  \\
\hspace{17mm}{\rm{subject}}\hspace{1mm}{\rm{to}}\hspace{3mm}\nabla g_{p,j}({\bf{u}}_k)^T ({\bf{u}} - {\bf{u}}_{k}) \leq -\delta_{g,j} \vspace{1mm} \\
\hspace{36mm}\forall j: g_{p,j}({\bf{u}}_k) \geq -\epsilon_j \vspace{1mm}\\
\hspace{36mm}{\bf{u}}^L \preceq {\bf{u}} \preceq {\bf{u}}^U
\end{array}.
\end{equation}

\noindent This is followed by replacing Condition (\ref{eq:suffcond4}) by the following line search in $K_k$:

\begin{equation}\label{eq:linesearch1}
\begin{array}{l}
\mathop {{\rm{minimize}}}\limits_{K_k, {\bf u}_{k+1}}\hspace{4mm}\phi ({\bf{u}}_{k+1}) \\
{\rm{subject}}\hspace{1mm}{\rm{to}}\hspace{3mm}  {\bf{u}}_{k+1}  = {\bf{u}}_k  + K_k \left( \bar {\bf{u}}^*_{k + 1} - {\bf{u}}_k \right) \\
\hspace{18mm}K_{k} \le \mathop{\min}\limits_{j = 1,...,n_g} \left[ \frac{\displaystyle{-g_{p,j}({\bf{u}}_{k } )}}{\displaystyle \sum\limits_{i = 1}^{n_u} {\kappa_{ji} | \bar u^*_{k+1,i} - u_{k,i} | }} \right] \\
\hspace{18mm}0 \leq K_k \leq 1 \\
\hspace{18mm}{\bf{u}}^L \preceq {\bf{u}}_{k+1} \preceq {\bf{u}}^U
\end{array}\hspace{3mm},
\end{equation}

\noindent which minimizes the cost along the projected direction while ensuring feasibility of all the uncertain plant constraints. We note that the projected direction need not be one of local descent for the cost, and that the values of $K_k$ that lead to a descent may consist of disjoint intervals due to multiple local minima. As such, this implementation is not guaranteed to converge to the local KKT point of the neighborhood in which it starts, but may be guaranteed to converge to some KKT point. To enforce this, we may simply fall back on using the original projection with the known cost:

\begin{equation}\label{eq:feasproj100}
\begin{array}{l}
\bar {\bf{u}}^*_{k+1} = {\rm{arg}} \mathop {{\rm{minimize}}}\limits_{{\bf{u}}}\hspace{4mm}\left\| {\bf{u}} - {\bf{u}}^*_{k+1} \right\|_2^2 \vspace{1mm}  \\
\hspace{17mm}{\rm{subject}}\hspace{1mm}{\rm{to}}\hspace{3mm}\nabla g_{p,j}({\bf{u}}_k)^T ({\bf{u}} - {\bf{u}}_{k}) \leq -\delta_{g,j} \vspace{1mm} \\
\hspace{36mm}\forall j: g_{p,j}({\bf{u}}_k) \geq -\epsilon_j \vspace{1mm}\\
\hspace{36mm}\nabla \phi({\bf{u}}_k)^T ({\bf{u}} - {\bf{u}}_{k}) \leq -\delta_\phi \vspace{1mm} \\
\hspace{36mm}{\bf{u}}^L \preceq {\bf{u}} \preceq {\bf{u}}^U
\end{array},
\end{equation}

\noindent instead of (\ref{eq:feasprojcostknown}) whenever (\ref{eq:linesearch1}) returns a step of $K_k \approx 0$ (i.e. when all steps along the projected direction either increase the cost or are infeasible), as doing this will ensure that we are able to find a local descent direction unless we are already arbitrarily close to a KKT point. Since (\ref{eq:linesearch1}) cannot continue to decrease the cost indefinitely, it follows that eventually we will be forced to return to the standard method and converge.

\subsubsection{Incorporating Some Known Constraint Functions}

Consider now the case where the cost is uncertain but some of the constraints are known:

\begin{equation}\label{eq:conknown}
\begin{array}{l}
\mathop {{\rm{minimize}}}\limits_{\bf{u}}\hspace{4mm}\phi_p ({\bf{u}}) \\
{\rm{subject}}\hspace{1mm}{\rm{to}}\hspace{3mm}{\bf{G}}_p({\bf{u}}) \preceq {\bf{0}} \\
\hspace{18mm}{\bf{G}}({\bf{u}}) \preceq {\bf{0}} 
\end{array}\hspace{3mm},
\end{equation}

\noindent where we have lumped the box constraints into the known constraint set.

Here, as before, we propose to remove the SCFO on the known parts from the projection, and simply carry out Projection (\ref{eq:feasproj}) without including the known constraints. $K_k$ is then chosen by the line search:

\begin{equation}\label{eq:linesearch2}
\begin{array}{l}
\mathop {{\rm{maximize}}}\limits_{K_k, {\bf u}_{k+1}}\hspace{4mm}K_k \\
{\rm{subject}}\hspace{1mm}{\rm{to}}\hspace{3mm}  {\bf{u}}_{k+1}  = {\bf{u}}_k  + K_k \left( \bar {\bf{u}}^*_{k + 1} - {\bf{u}}_k \right) \\
\hspace{18mm}K_{k} \le \mathop{\min}\limits_{j = 1,...,\bar n_g} \left[ \frac{\displaystyle{-g_{p,j}({\bf{u}}_{k } )}}{\displaystyle \sum\limits_{i = 1}^{n_u} {\kappa_{ji} | \bar u^*_{k+1,i} - u_{k,i} | }} \right] \\
\hspace{18mm}K_k < -2\displaystyle\frac{\nabla \phi_p({\bf{u}}_k)^T (\bar{\bf{u}}^*_{k+1} - {\bf{u}}_{k})}{(\bar{\bf{u}}^*_{k+1} - {\bf{u}}_{k})^T \overline {\bf{Q}}_\phi (\bar{\bf{u}}^*_{k+1} - {\bf{u}}_{k})} \\
\hspace{18mm}K_k \leq 1 \\
\hspace{18mm}{\bf{G}}({\bf{u}}_{k+1}) \preceq {\bf{0}}
\end{array}\hspace{3mm},
\end{equation}

\noindent with $\bar n_g$ used to denote the number of uncertain constraints ($n_g$ denoting the \emph{total} number of known and uncertain ones).

It may be that (\ref{eq:linesearch2}) returns $K_k = 0$ if the projected direction is infeasible for all $K_k \in (0,1]$ for the combination of the known constraints (a necessary condition for this is that at least one of them be active). In this case, once again, we may fall back on the original methodology and incorporate these particular constraints into the projection:

\begin{equation}\label{eq:feasprojconknown}
\begin{array}{l}
\bar {\bf{u}}^*_{k+1} = {\rm{arg}} \mathop {{\rm{minimize}}}\limits_{{\bf{u}}}\hspace{4mm}\left\| {\bf{u}} - {\bf{u}}^*_{k+1} \right\|_2^2 \vspace{1mm}  \\
\hspace{17mm}{\rm{subject}}\hspace{1mm}{\rm{to}}\hspace{3mm}\nabla g_{p,j}({\bf{u}}_k)^T ({\bf{u}} - {\bf{u}}_{k}) \leq -\delta_{g,j} \vspace{1mm} \\
\hspace{36mm}\forall j \in [1,\bar n_g]: g_{p,j}({\bf{u}}_k) \geq -\epsilon_j \vspace{1mm}\\
\hspace{36mm}\nabla g_{j}({\bf{u}}_k)^T ({\bf{u}} - {\bf{u}}_{k}) \leq -\delta_{j} \vspace{1mm} \\
\hspace{36mm}\forall j \in [\bar n_g+1,  n_g]: g_{j}({\bf{u}}_k) = 0 \vspace{1mm}\\
\hspace{36mm}\nabla \phi_{p}({\bf{u}}_k)^T ({\bf{u}} - {\bf{u}}_{k}) \leq -\delta_\phi \vspace{1mm} \\
\hspace{36mm}{\bf{u}}^L \preceq {\bf{u}} \preceq {\bf{u}}^U
\end{array}.
\end{equation}

Note that we are not compelled to use (\ref{eq:feasprojconknown}) by default -- that is, without trying (\ref{eq:feasproj}) first -- as there may be instances when a known constraint is active but (\ref{eq:feasproj}) yields a direction for which it may remain feasible. We also note, however, that (\ref{eq:feasprojconknown}) is sufficient to enforce the original guarantees, as its not having a solution may be seen as proof of KKT convergence (to a point that is potentially defined by a combination of known and uncertain constraints).

\subsubsection{Incorporating a Known Cost and Some Known Constraint Functions}

Finally, we investigate the case where both the cost and some constraints are known:

\begin{equation}\label{eq:concostknown}
\begin{array}{l}
\mathop {{\rm{minimize}}}\limits_{\bf{u}}\hspace{4mm}\phi ({\bf{u}}) \\
{\rm{subject}}\hspace{1mm}{\rm{to}}\hspace{3mm}{\bf{G}}_p({\bf{u}}) \preceq {\bf{0}} \\
\hspace{18mm}{\bf{G}}({\bf{u}}) \preceq {\bf{0}}
\end{array}\hspace{3mm},
\end{equation}

\noindent for which the following line search is proposed:

\begin{equation}\label{eq:linesearch3}
\begin{array}{l}
\mathop {{\rm{minimize}}}\limits_{K_k, {\bf u}_{k+1}}\hspace{4mm}\phi ({\bf{u}}_{k+1}) \\
{\rm{subject}}\hspace{1mm}{\rm{to}}\hspace{3mm}  {\bf{u}}_{k+1}  = {\bf{u}}_k  + K_k \left( \bar {\bf{u}}^*_{k + 1} - {\bf{u}}_k \right) \\
\hspace{18mm}K_{k} \le \mathop{\min}\limits_{j = 1,...,\bar n_g} \left[ \frac{\displaystyle{-g_{p,j}({\bf{u}}_{k } )}}{\displaystyle \sum\limits_{i = 1}^{n_u} {\kappa_{ji} | \bar u^*_{k+1,i} - u_{k,i} | }} \right] \\
\hspace{18mm}0 \leq K_k \leq 1 \\
\hspace{18mm}{\bf{G}}({\bf{u}}_{k+1}) \preceq {\bf{0}}
\end{array}\hspace{3mm},
\end{equation}

\noindent and may be done after sequentially attempting Projections (\ref{eq:feasprojcostknown}) $\rightarrow$ (\ref{eq:feasproj100}) $\rightarrow$ (\ref{eq:feasproj200}), with the latter given as:

\begin{equation}\label{eq:feasproj200}
\begin{array}{l}
\bar {\bf{u}}^*_{k+1} = {\rm{arg}} \mathop {{\rm{minimize}}}\limits_{{\bf{u}}}\hspace{4mm}\left\| {\bf{u}} - {\bf{u}}^*_{k+1} \right\|_2^2 \vspace{1mm}  \\
\hspace{17mm}{\rm{subject}}\hspace{1mm}{\rm{to}}\hspace{3mm}\nabla g_{p,j}({\bf{u}}_k)^T ({\bf{u}} - {\bf{u}}_{k}) \leq -\delta_{g,j} \vspace{1mm} \\
\hspace{36mm}\forall j \in [1,\bar n_g]: g_{p,j}({\bf{u}}_k) \geq -\epsilon_j \vspace{1mm}\\
\hspace{36mm}\nabla g_{j}({\bf{u}}_k)^T ({\bf{u}} - {\bf{u}}_{k}) \leq -\delta_{j} \vspace{1mm} \\
\hspace{36mm}\forall j \in [\bar n_g+1, n_g]: g_{j}({\bf{u}}_k) = 0 \vspace{1mm}\\
\hspace{36mm}\nabla \phi({\bf{u}}_k)^T ({\bf{u}} - {\bf{u}}_{k}) \leq -\delta_\phi \vspace{1mm} \\
\hspace{36mm}{\bf{u}}^L \preceq {\bf{u}} \preceq {\bf{u}}^U
\end{array}.
\end{equation}

As a final note to this discussion, we may remark that, if the cost is known and strictly convex (e.g. positive-definite quadratic), then:

\begin{equation}\label{eq:convexity1}
\phi (\bar {\bf{u}}^*_{k+1}) < \phi ({\bf{u}}_{k}) \Rightarrow \phi ({\bf{u}}_{k+1}) < \phi ({\bf{u}}_{k}),\;\; \forall K_k \in (0,1],
\end{equation}

\noindent and that, likewise, if a constraint is known and convex:

\begin{equation}\label{eq:convexity2}
g (\bar {\bf{u}}^*_{k+1}) \leq 0 \Rightarrow g ({\bf{u}}_{k+1}) \leq 0,\;\; \forall K_k \in [0,1].
\end{equation}

\noindent Both of these statements are easily proven by Jensen's inequality.

(\ref{eq:convexity2}) essentially explains why we did not put the box constraints through the same analysis as we did general known constraints -- as the RTO algorithm and the projected target both satisfy them by default and as they are convex, any choice of the filter gain between 0 and 1 is sufficient to keep them satisfied.

\vspace{2mm}
\noindent{\it{Examples with Known Elements}} 
\vspace{2mm}

\begin{figure}
\begin{center}
\includegraphics[width=8cm]{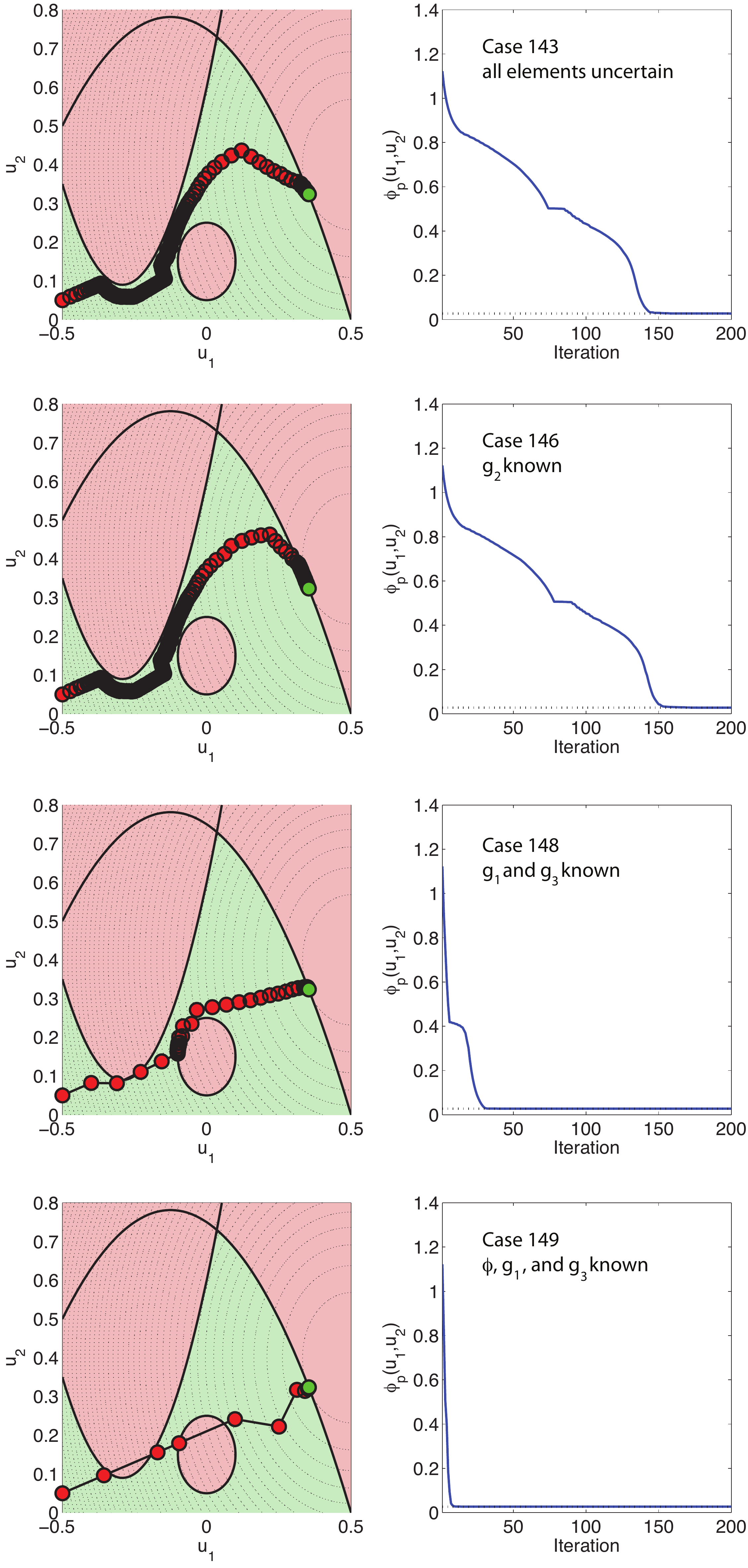}    
\caption{Performance of the ideal-target algorithm for Problem A when partial knowledge is assumed.}
\label{fig:P2ex15}
\end{center}
\end{figure}

To illustrate the benefits of incorporating known elements, we consider Problem A with the ideal target scheme over 200 iterations. To highlight the benefits of a known cost, we also use a more conservative $\overline {\bf Q}_\phi = 20{\bf I}$ instead of the cost function Hessian ($2 {\bf I}$). Multiple realizations, with knowledge assumed in different combinations of cost/constraints, are considered and reported in Table \ref{tab:certain} and Figure \ref{fig:P2ex15}.

From these results, we first see that knowledge does not automatically improve performance and may, in some cases, even make it worse by changing the path of the algorithm (e.g. Case 146 as compared to Case 143). This is, however, mostly true when the known parts are not important factors in limiting convergence speed, and it can be seen in the example that making either one or both of the constraints $g_{p,1}$ and $g_{p,3}$ known leads to an immense improvement in performance. Additionally, while adding a known cost does not make much of a difference when all of the constraints are uncertain (compare Case 143 and Case 144) since the algorithm is anyway forced to take small steps to remain feasible, doing so when knowledge of both $g_{p,1}$ and $g_{p,3}$ is assumed leads to even greater benefits (compare Cases 148 and 149). For the latter, the algorithm is able to virtually jump over these constraints without having to maneuver around them (bottom of Figure \ref{fig:P2ex15}).

\begin{table}
\caption{Performance of the ideal-target scheme for Problem A when different parts of the problem are assumed to be known.}
\begin{center}
\footnotesize{
\begin{tabular}{|p{.6cm}|p{1.5cm}|p{1cm}|}
\hline
Case & Known & $L$   \\ \hline
143 & None & 76.86 \\ \hline
144 & $\phi$ & 76.20 \\ \hline 
145 & $g_{1}$ & 26.99 \\ \hline 
146 & $g_{2}$ & 81.51 \\ \hline
147 & $g_{3}$ & 65.40 \\ \hline
148 & $g_{1}$, $g_{3}$ & 9.78 \\ \hline
149 & $\phi$, $g_{1}$, $g_{3}$ & 3.00 \\ \hline
\end{tabular}
}
\end{center}
\label{tab:certain}
\end{table}

\section{Summary, Reflections, and an Outlook on the Future}

This paper has extended on the work of its predecessor \citep{Bunin:12b}, where the sufficient conditions for feasibility and optimality for the general RTO algorithm were presented but were ultimately deemed inapplicable as they relied on knowledge that is unavailable in a practical setting. Here, we have systematically treated these ``knowledge gaps'' one by one -- proposing robust implementations of the SCFO for the cases where the true plant gradients and constraint values were replaced by bounded estimates, and then proposing some guidelines by which the Lipschitz and quadratic upper bounds could be chosen and then refined so as to reduce their innate conservatism. Additionally, we have demonstrated that incorporating partial analytical knowledge into the uncertain RTO problem, as well as allowing the violation of certain constraints, is both straightforward and can lead to significant performance benefits with regard to convergence speed.

It is interesting to note that the different implementation issues of the SCFO are not issues inherent just to the SCFO but are, in fact, standard engineering problems. As such, progress in these directions will only aid in RTO performance. The gradient estimation problem of obtaining good derivative estimates from discrete data, for example, is a research topic where any advances would serve to tighten the gradient uncertainty bounds and to thereby, by Theorem 1, reduce the size of the KKT neighborhood to which an SCFO-supplemented algorithm would converge. Likewise, improvements in sensor technology, modeling, and output estimation methods are all likely to lead to better bounds on constraint value estimates (the $\overline g_{p,j} ({\bf u}_k)$ in Section 3) -- thereby diminishing the decrease in convergence speed that results when these bounds are too conservative. Finally, modeling and enhanced process knowledge, together with robust global optimization methods to solve problems like (\ref{eq:parlip}), will all contribute to better offline determination of the Lipschitz constants and quadratic upper bounds. We also note that this latter is a completely new topic of its own, as the authors are not familiar with any research that has attempted to systematically find these (very useful) constants in any application, and as such the contents of Section 4 only represent an initial attempt at some reasonable approaches.

A number of other important questions remain unanswered and should eventually be the subject of a future study. As hinted at in the end of the previous paper and as corroborated to an extent by the numerical studies of Section 2 here, it is unclear when the application of the SCFO is beneficial, as a projection of a target given by the RTO algorithm based only on local information can be detrimental to convergence speed, especially when the RTO target is very good. As such, while the SCFO make RTO more \emph{reliable}, they may not necessarily make it more performant (though in certain cases with difficult constraints like in Problem A, or in any problem with a poor RTO algorithm, they can). Furthermore, it is unclear as to how such performance losses may scale with an increasing number of input variables. Both of these questions require further study.

Three crucial implementation issues that have not been addressed here and that immediately come to mind are the following:

\begin{itemize}
\item The handling of noise or errors in the \emph{inputs}, which here have been assumed to be manipulated perfectly. In some applications, the inputs themselves may be subject to uncertainty -- returning to the example of the CSTR in Section 4, consider the case where disturbances in the piping do not allow for perfect control, or measurement, of $T_j$.
\item The minimal-excitation requirement that consecutive inputs, ${\bf u}_k$ and ${\bf u}_{k+1}$, be far enough apart so as to make tasks like gradient or parameter estimation possible. Obviously, such a condition and monotonic convergence to a KKT point (as promised by the SCFO) are in conflict since convergence implies ${\bf u}_k = {\bf u}_{k+1}$.
\item The topic of process degradation, as (almost) all of the theory proposed in these two papers has only focused on cases where the RTO problem remains the same for all iterations and where $\phi_p ({\bf{u}})$ and ${\bf{G}}_p({\bf{u}})$ do not change. In practice, equipment may degrade, environmental conditions may change, and various disturbances may enter to gradually (or not so gradually) change the RTO problem over the course of operation.
\end{itemize}

\noindent Proper, rigorous treatment of these problems could constitute a third publication in this series. While challenging, we are optimistic that the flexibility of the SCFO framework, as evidenced in this paper, will allow room for these scenarios as well.

Finally, it would be interesting to see how SCFO-supplemented RTO schemes perform in real experimental problems! While a rudimentary application of the SCFO to an unconstrained problem of iterative controller tuning has already given positive results \citep{Bunin:13}, it is of great interest to see the results for all sorts of RTO problems, and several other applications are already in planning, together with an open-source solver that incorporates all of the theory discussed in this paper \citep{SCFOug}.

\bibliography{feasconv}             
                                                   
\end{document}